\documentclass[a4paper,12pt]{article}

\usepackage{amsfonts,amsmath,amsopn,amsthm}
\usepackage{amssymb,mathrsfs}
\usepackage{hyperref}

\makeatletter
\@addtoreset{equation}{section}

\makeatother

\newtheorem{lemma}{Lemma}[section]
\newtheorem{corollary}[lemma]{Corollary}
\newtheorem{theorem}[lemma]{Theorem}

\theoremstyle{remark}
\newtheorem{massumption}[lemma]{Main assumption}
\newtheorem{sassumption}[lemma]{Secondary assumption}
\newtheorem{remark}[lemma]{Remark}
\newtheorem*{claim}{Claim}

\newcommand{\bC}{\mathbf{C}}
\newcommand{\bR}{\mathbf{R}}
\newcommand{\bZ}{\mathbf{Z}}
\newcommand{\cA}{\mathscr{A}}
\newcommand{\cE}{\mathscr{E}}
\newcommand{\cF}{\mathscr{F}}
\newcommand{\cG}{\mathscr{G}}
\newcommand{\cM}{\mathscr{M}}
\newcommand{\cO}{\mathscr{O}}
\newcommand{\fu}{\mathfrak{u}}

\newcommand{\iA}{{^i\!A}}
\newcommand{\iB}{{^i\!B}}
\newcommand{\lt}{\mathfrak{t}}
\newcommand{\cL}{\mathscr{L}}

\newcommand{\Lie}{\mathrm{Lie}}
\newcommand{\zb}{\bar{z}}
\newcommand{\db}{\bar{\partial}}

\newcommand{\spq}{/\!\!/}
\renewcommand{\geq}{\geqslant}
\renewcommand{\leq}{\leqslant}

\DeclareMathOperator{\im}{Im}
\DeclareMathOperator{\re}{Re}
\DeclareMathOperator{\ad}{ad}
\DeclareMathOperator{\tr}{tr}
\DeclareMathOperator{\rank}{rank}
\DeclareMathOperator{\diag}{diag}
\DeclareMathOperator{\pdeg}{p-deg}
\DeclareMathOperator{\End}{End}
\DeclareMathOperator{\coker}{coker}
\DeclareMathOperator{\scal}{scal}

\author{Olivier Biquard and Philip Boalch\footnote{Both authors are members of EDGE,
Research Training Network HPRN-CT-2000-00101, supported by the
European Human Potential Programme.}}
\date{}
\title{Wild nonabelian Hodge theory on curves}

\begin{document}
\maketitle

\begin{abstract}
On a complex curve, we establish a correspondence between integrable
connections with irregular singularities, and Higgs bundles such that
the Higgs field is meromorphic with poles of any order.

The moduli spaces of these objects are obtained by fixing at each
singularity the polar part of the connection, which amounts to fixing
a coadjoint orbit of the group $GL_r(\bC[z]/z^n)$. We prove
that they carry hyperK{\"a}hler metrics, which are complete when the
residues of the connection are semisimple. 
\end{abstract}

\section*{Introduction}

On a compact K{\"a}hler manifold, there is a well-known correspondence,
named nonabelian Hodge theory, and established by Simpson and
Corlette (see \cite{Si92}), between representations of the
fundamental group (or integrable connections) and Higgs bundles.
In the case of a curve, the correspondence is due to Hitchin
\cite{Hi87} and Donaldson \cite{Do87}.
This correspondence has been extended to the case when the
objects on both sides have logarithmic singularities, at least in the
case when the singular locus is a smooth divisor, see \cite{Si90} on
curves and \cite{Bi97} in higher dimensions.

In this article, we extend the correspondence to the irregular case,
on a curve. This means that we now look, on one side, at integrable
connections with irregular singularities, like
\begin{equation}\label{int-D}
d+A_n \frac{dz}{z^n} + \cdots + A_1 \frac{dz}{z}
\end{equation}
with $n>1$, and, on the other side, Higgs bundles $(\cE,\theta)$ with
the Higgs field having a polar part
\begin{equation}\label{int-H}
\theta = T_n \frac{dz}{z^n} + \cdots + T_1 \frac{dz}{z}.
\end{equation}
In this paper, we will restrict to the case when one can choose
\begin{equation}\label{ss}
A_2,\dots,A_n \textrm{ semisimple.}
\end{equation}

Sabbah \cite{Sa99} has constructed a harmonic metric for integrable
connections with irregular singularities: this is a part of the
correspondence. We construct the whole correspondence:
\begin{theorem}\label{th:cor}
Under the assumption (\ref{ss}),
there is a 1-1 correspondence between stable (parabolic) integrable
connections with irregular singularities, and stable parabolic Higgs
bundles with singularities like in (\ref{int-H}).
The correspondence between the singularities (\ref{int-D}) and
(\ref{int-H}) at the punctures is given by 
$$ T_i=\frac{1}{2} A_i, \quad i\geq 2. $$
For $i=1$ there is the same permutation between eigenvalues of $A_1$,
$T_1$, and the parabolic weights as in \cite{Si90}: more precisely if
on the connection side we have the eigenvalues $\mu_i$ of $A_1$ and
the parabolic weights $\beta_i$, and on the Higgs bundle side the
eigenvalues $\lambda_i$ of $T_1$ and the parabolic weights $\alpha_i$, then
$$ \alpha_i=\re\mu_i-[\re\mu_i],\quad
\lambda_i=\frac{\mu_i-\beta_i}{2}, $$
where the $[\cdot]$ denotes the integer part.
The nilpotent parts of $A_1$ and $T_1$ are the same.
\end{theorem}
See section \ref{s_complex} for details about stability. When the
weights of the local system vanish ($\re\lambda_i=0)$, stability for integrable
connections reduces to irreducibility, see remark \ref{rem:weights}.
The proof also gives precise information about the asymptotics of the
harmonic metric.

In \cite{Bo99} the second author studied the symplectic geometry of
the moduli spaces of integrable connections with irregular
singularities in the case when $A_n$ is regular:
more precisely, the moduli space obtained by fixing the coadjoint orbit of
the group $GL_r(\bC[z]/z^n)$ defined by (\ref{int-D}) is
complex symplectic.
In more concrete terms, this means considering the moduli space of
connections which near the singularity can be transformed into 
$$
d+A_n \frac{dz}{z^n} + \cdots + A_1 \frac{dz}{z}+\textrm{holomorphic terms.}
$$
In the case where $A_n$ is regular semisimple,
this amounts to fixing the formal type of the connection.
Now we will show that this moduli space carries more structure, namely:
\begin{theorem}\label{th:moduli}
The moduli space of integrable connections, with fixed polar part
(\ref{int-D}) at each singularity, is hyperK{\"a}hler.
If the moduli space is smooth and the residues $A_1$ at each singular
point are semisimple, then the metric is complete. 
\end{theorem}
The same result remains true if we add some compatible parabolic
structure at each singularity. We remark that a generic choice of
parabolic structure leads to smooth moduli spaces (no semistable
points).

Also, at least generically, one can explicitly describe the underlying
complex symplectic manifold:
\begin{theorem}\label{th:generic}
For generic eigenvalues of the residues, all integrable connections
are stable, and the moduli space can be identified, from the
complex symplectic viewpoint, 
with the finite dimensional quasi-Hamiltonian quotient of \cite{saqh}.
\end{theorem}
For example over the projective line the moduli space contains
a dense open subset, parameterising connections on trivial holomorphic 
bundles,
which may be described as a complex symplectic quotient
of finite dimensional coadjoint orbits. However in general there are stable
connections on nontrivial holomorphic bundles 
and  the quasi-Hamiltonian quotient
incorporates these points as well.
An example will be given in section \ref{s_DR}.

Note that Martinet-Ramis \cite{MarRam91} have constructed a ``wild
fundamental group'', so that connections with irregular singularities
can be interpreted as finite dimensional representations of this
group. From this point of view, theorem \ref{th:cor} really generalizes the
earlier correspondences between representations of the fundamental
group of the curve (or the punctured curve) and Higgs bundles.

Starting on the other side we remark that
moduli spaces of meromorphic Higgs bundles have previously been
studied algebraically by Bottacin \cite{Bottacin} and Markman \cite{Markman}, 
who have shown they are algebraic completely integrable systems.

Also the study of Higgs bundles with such irregular singularities has
a physical interest, since some of them arise from the Nahm transform of
periodic monopoles, see \cite{CheKap01}. The corresponding hyperK{\"a}hler
metrics have been studied in \cite{CheKap}.

It is not clear at all how to extend the correspondence to the case
where $A_2,\dots,A_n$ have a nilpotent part. For example, in
the case of an order 2 pole, if we suppose that $T_2$ in (\ref{int-H})
is nilpotent, then the eigenvalues of $\theta$ have only a simple
pole, so the Higgs bundle satisfies the ``tameness'' condition of
Simpson \cite{Si90}, and we are back in the simple pole situation.
In particular in that case, we cannot construct new metrics:
actually, by a meromorphic gauge transformation, the Higgs field can
be transformed into a Higgs field with only a simple pole. The same
kind of considerations holds on the integrable connection side.

One of the main features of connections with irregular singularities
is that formal equivalence does not come from holomorphic equivalence,
resulting in the well-known Stokes phenomenon.
Sabbah studies carefully this Stokes phenomenon around the puncture in
order to construct a sufficiently good initial metric to which he can
apply Simpson's existence theorem for harmonic metrics \cite{Si90}.

Our method is different: we use a weighted Sobolev space approach, which
enables us to forget the difficult structure of irregular
singularities, at the expense of developing some analysis to
handle the operators with highly singular coefficients that we encounter.
In particular, we strengthen the Fr{\'e}chet symplectic quotient of \cite{Bo99} 
into a hyperK{\"a}hler quotient.

In section \ref{s_model}, we develop the local models which are
a guide for the behavior of the correspondence, and then define the
admissible deformations in suitable Sobolev spaces in section
\ref{s_deformations}. Next we study the local analysis needed on a
disk in sections \ref{s_gauge} and \ref{s_analysis}.
This enables us to construct the $C^\infty$-moduli space of solutions
of Hitchin's selfduality equations in section \ref{s_moduli}, and
prove that it is hyperK{\"a}hler. It remains to identify this moduli space
with the ``DeRham moduli space'' of integrable connections and the
``Dolbeault moduli space'' of Higgs bundles: these moduli spaces are
studied in sections \ref{s_Dol} and \ref{s_DR}, the correspondence is
stated in section \ref{s_complex} and proven in section \ref{s_proof}.

\noindent\emph{Acknowledgments.}
We thank Claude Sabbah for discussions and interest and 
Andrey Bolibruch for useful comments.

\section{Local model}\label{s_model}

Look at a rank $r$ holomorphic bundle $\cF$ in the unit disk,
trivialized in a basis $(\tau_i)$ of holomorphic sections, and
consider the holomorphic connection
\begin{equation}\label{m_D}
D = d + A_n \frac{dz}{z^n} + \cdots + A_1 \frac{dz}{z},
\end{equation}
where the $A_i$ are constant matrices; then $D$ is an
integrable connection, that is, the curvature $F_D=D^2$ vanishes.

The work in this paper will be valid under the following assumption.
\begin{massumption}\label{ass1}
The $(A_i)_{i\geq 2}$ are diagonal matrices.
\end{massumption}

Notice that this implies that the $A_i$ commute.
Actually, it would be sufficient to suppose that the $A_i$ are semisimple.
Indeed, under this weaker hypothesis, one can do the following:
observe that a holomorphic transformation $e^u$ transforms $D$ into
$$
d + e^{\ad u}\big(A_n \frac{dz}{z^n} + \cdots + A_1 \frac{dz}{z})
  - (de^u)e^{-u} ;  
$$
therefore, if $A_n$ is semisimple, we can find $u$ holomorphic such
that $u(0)=0$ and
$$ e^u(D)=d+ A_n \frac{dz}{z^n}
           + A'_{n-1} \frac{dz}{z^{n-1}} + \cdots + A'_1 \frac{dz}{z}
           + \textrm{ holomorphic terms,}
$$
and $[A_n,A'_i]=0$. This means that, up to holomorphic terms, one can
suppose that the $(A_i)_{i\leq n-1}$ commute with $A_n$.
Iterating this procedure, we see that, up to holomorphic
terms, $D$ is equivalent to a connection (\ref{m_D}) where the
matrices $A_i$ commute.

If none of the $(A_i)_{i\geq 2}$ are regular (more precisely if the
stabilizer of $(A_i)_{i\geq 2}$ is not a
Cartan subalgebra) then the residue $A_1$ may have a nilpotent part. The
phenomenon induced by such a nilpotent part in the residue of the
connection $D$ is complicated, but has been carefully analyzed in
\cite{Si90,Bi97}.
As the purpose of this work is mainly to investigate the new
phenomenon induced by the presence of the higher order terms
$(A_i)_{i\geq 2}$, and since there is no interaction with the
nilpotent part of the residue (actually the two behaviors, coming from
the higher order poles on one side, from the nilpotent part of the
residue on the other side, can basically be superposed), we will make
another assumption for the clarity of exposition. The
results remain true without the assumption, but the precise proof
would require incorporating the analysis in \cite{Bi97} into that
developed here: this is not difficult (see remark
\ref{r_nil}), but would lead to very heavy definitions of functional
spaces, so we prefer to give a proof of the results only under the
following hypothesis.
\begin{sassumption}\label{ass2}
The matrix $A_1$ is also a diagonal matrix.
\end{sassumption}

Finally we will suppose that $\cF$ comes with a parabolic structure with
weights $\beta_i \in [0,1[$, meaning basically that we
have on the bundle $\cF$ a hermitian metric
\begin{equation}
h=\left(\begin{array}{ccc}|z|^{2\beta_1} & & \\ & \ddots & \\ & & |z|^{2\beta_r}
	\end{array}
\right) ;
\end{equation}
the fiber $\cF_0$ at the origin is filtered by
$\cF_\beta=\{s(0),|s(z)|=O(|z|^\beta)\}$.
In the orthonormal basis $(\tau_i/|z|^{\beta_i})$, we get the formula
\begin{equation}\label{f_D}
D = d + \sum_1^n A_i \frac{dz}{z^i} - \beta \frac{dr}{r},
\end{equation}
where $\beta$ is the diagonal matrix with coefficients
$\beta_1,\ldots,\beta_r$.

Recall that in general, we have a decomposition of $D$ into a unitary
part and a selfadjoint part,
$$ D = D^+ + \phi , $$
and we can define new operators
\begin{eqnarray*}
D'' &=& (D^+)^{0,1} + \phi^{1,0}, \\
D'  &=& (D^+)^{1,0} + \phi^{0,1}.
\end{eqnarray*}
The operator
\begin{equation}\label{f_higgs}
D'' = \db^E + \theta
\end{equation}
is a candidate to define a Higgs bundle structure, and this is the
case if the \emph{pseudo-curvature} 
\begin{equation}
G_D = -2(D'')^2
\end{equation}
vanishes. In the case of a Riemann surface, the equation reduces to
$\db^E\theta=0$.

In our case, we get, still in the orthonormal basis $(\tau_i/|z|^{\beta_i})$,
\begin{eqnarray*}
D^+ &=&
 d + \frac{1}{2} \sum_1^n A_i \frac{dz}{z^i} - A_i^* \frac{d\zb}{\zb^i}, \\
\phi&=&
 \frac{1}{2} \sum_1^n A_i \frac{dz}{z^i} + A_i^* \frac{d\zb}{\zb^i}
         - \beta \frac{dr}{r},\\
\db^E&=& \db -\frac{1}{2} \sum_1^n A_i^* \frac{d\zb}{\zb^i}, \\
\theta&=& \frac{1}{2}\sum_1^n A_i \frac{dz}{z^i}
         - \frac{\beta}{2} \frac{dz}{z} .
\end{eqnarray*}
It is clear that $G_D=0$, so that we have a solution of
Hitchin's selfduality equations.

Formulas become simpler if we replace the orthonormal basis
$(\tau_i/|z|^{\beta_i})$ by the other orthonormal basis $(e_i)$ given
on the punctured disk by
$$
e_i = |z|^{-\im \mu_i} \exp\bigg(
  \sum_2^n -\frac{A_i^*}{2(i-1)\zb^{i-1}} + \frac{A_i}{2(i-1)z^{i-1}}
          \bigg) \frac{\tau_i}{|z|^{\beta_i}},
$$
where the $\mu_i$ are the eigenvalues of the residue of $D$ (the
diagonal coefficients of $A_1$). 
Indeed, in the basis $(e_i)$, the previous formulas become 
\begin{equation}\label{f_modDR}
D^+ = d+\re(A_1)i d\theta, \quad
\phi=\frac{1}{2} \sum_1^n A_i \frac{dz}{z^i} + A_i^* \frac{d\zb}{\zb^i}
         - \beta \frac{dr}{r}
\end{equation}
and
\begin{equation}\label{f_modDol}
\db^E=\db-\frac{1}{2}\re(A_1)\frac{d\zb}{\zb}, \quad 
\theta= \frac{1}{2}\sum_1^n A_i \frac{dz}{z^i}
         - \frac{\beta}{2} \frac{dz}{z} .
\end{equation}
This orthonormal basis $(e_i)$ defines a hermitian extension $E$ of
the bundle over the puncture.

Now look at the holomorphic bundle induced by $\db^E$ on the punctured
disk. A possible choice of a basis of holomorphic sections is
\begin{equation}\label{f_sigma}
\sigma_i = |z|^{\re\mu_i-[\re\mu_i]} e_i.
\end{equation}
We see that $|\sigma_i| = |z|^{\alpha_i}$,
with
\begin{equation}\label{f_alpha}
\alpha_i = \re\mu_i-[\re\mu_i].
\end{equation}
This choice of $\sigma_i$ is the only possible choice for which
$0\leq\alpha_i<1$. The sections $(\sigma_i)$ define an extension $\cE$ of
the holomorphic bundle over the puncture, and the behavior of the
metric means that this extension carries a parabolic structure with
weights $\alpha_i$. Finally, in this basis the Higgs field is still
given by
\begin{equation}\label{f_theta}
\theta = \frac{1}{2}\sum_1^n A_i \frac{dz}{z^i}
        - \frac{\beta}{2} \frac{dz}{z}.
\end{equation}
In particular, the eigenvalues of the residue of the Higgs field are
\begin{equation}\label{f_lambda}
\lambda_i = \frac{\mu_i-\beta_i}{2}.
\end{equation}
The formulas (\ref{f_alpha}) and (\ref{f_lambda}) give the same
relations between parabolic weights and eigenvalues of the residue on
both sides as in the case of regular singularities. This basically
means that the behavior described by Simpson in the case of regular
singularities still occurs here in the background of the behavior of
the solutions in presence of irregular singularities.

\section{Deformations}\label{s_deformations}

We continue to consider connections in a disk, using the same
notations as in section \ref{s_model}.

We want to construct a space $\cA$ of admissible connections on $E$,
with the same kind of singularity as $D$ at the puncture. In order
to be able to do some analysis, we need to define Sobolev
spaces. First define a weighted $L^2$ space (using the function $r=|z|$)
$$
L^2_\delta = \{ f, \frac{f}{r^{\delta+1}}\in L^2 \}.
$$
The convention for the weight is chosen so that the function $r^x \in
L^2_\delta$ if and only if $x>\delta$.

Now we want to define Sobolev spaces for sections $f$ of $E$ or of
associated bundles, mainly $\End(E)$. Let us restrict to this case: we
have to decompose $\End(E)$ under the action of the $A_i$.

A simple case is the \emph{regular case} in which the stabilizer of $A_n$ is
the same as the stabilizer of all the matrices $A_1,\dots,A_n$. Then
we decompose $\End(E)$ as
\begin{equation}
\End(E)=\End(E)_n \oplus \End(E)_0,\quad
\left\{\begin{array}{l}
\End(E)_0=\ker\ad(A_n),\\ \End(E)_n=(\ker\ad(A_n))^\perp
       \end{array}\right..
\end{equation}
For example, if $A_n$ is regular, then $\End(E)_0$ consists of diagonal
matrices.

In the non regular case, we need a more subtle decomposition,
\begin{equation}
\End(E)=\sum_0^n \End(E)_i,
\end{equation}
defined by induction by
\begin{equation}
\End(E)_0=\cap \ker\ad(A_i),\quad
\End(E)_i=\End(E)_{i-1}^\perp \cap \big(\cap_{j>i} \ker\ad(A_j)\big).
\end{equation}
We will therefore decompose a section $f$ of $\End(E)$ as
$$
f=f_0+\cdots+f_n,
$$
where the indices mean that the highest order pole term acting on $f_i$ is
$A_i dz/z^i$.

Now we can define Sobolev spaces with $k$ derivatives in $L^2$,
\begin{equation}
L^{k,2}_\delta =
\{ f, \frac{\nabla^j f_i}{r^{i(k-j)}} \in L^2_\delta \textrm{ for }0\leq j\leq
k, 0\leq i\leq n \};
\end{equation}
in the whole paper, $\nabla=\nabla^{D^+}$ is the covariant derivative
associated to the \emph{unitary} connection $D^+$.

In this problem it is natural to look at deformations of $D$ such that
the curvature remains $O(r^{-2+\delta})$, that is slightly better than
$L^1$. This motivates the following definition of the space $\cA$ of
admissible deformations of $D$:
\begin{equation}
\cA = \{ D+a, a\in L^{1,2}_{-2+\delta}(\Omega^1\otimes\End E) \},
\end{equation}
and of the gauge group,
\begin{equation}\label{d_gauge}
\cG = \{ g\in U(E), Dg g^{-1}\in L^{1,2}_{-2+\delta} \}.
\end{equation}
The following lemma says that we have defined good objects for gauge theory.
\begin{lemma}\label{l_gt}
The connections in $\cA$ have their curvature in $L^2_{-2+\delta}$.
Moreover, $\cG$ is a Lie group, with Lie algebra
\begin{equation}
\Lie(\cG)= L^{2,2}_{-2+\delta}(\fu(E)),
\end{equation}
and it acts smoothly on $\cA$.
\end{lemma}
The proof of the lemma involves some nonlinear analysis, which we will
develop in section \ref{s_gauge}.

\begin{remark}\label{r_nil}
In the case $A_1$ has a nilpotent part, the analysis has to be refined
as in \cite{Bi97} in logarithmic scales in order to allow the
curvature to be $O(r^{-2}|\ln r|^{-2-\delta})$. The analysis which
will be developed here shall continue to be valid in this case, for
components in $\End(E)_n\oplus\cdots\oplus\End(E)_2$, and the tools in
\cite{Bi97} handle $\End(E)_1\oplus\End(E)_0$, where the action of
the irregular part is not seen. This is the basic reason why the results
in this paper continue to be true in this case also.
\end{remark}

\section{Gauge theory}\label{s_gauge}
In this section, we give the tools to handle the nonlinearity in the
gauge equations for connections with irregular singularities. 
We first need to define weighted $L^p$-spaces, for sections $f$ of $\End E$, by
$$
L^p_\delta = \{ f, \frac{f}{r^{\delta+2/p}} \in L^p \},\quad
L^{k,p}_\delta = \{f, \frac{\nabla^jf_i}{r^{i(k-j)}}\in L^p_\delta \}.
$$
It is convenient to introduce spaces with a weighted condition on
$f_0$ also,
$$
\hat{L}^{k,p}_\delta = \{f\in L^{k,p}_\delta,
\frac{\nabla^jf_0}{r^{k-j}}\in L^p_\delta \}.
$$
Again, the weight is chosen so that $r^x\in L^p_\delta$ if and only if
$x>\delta$. This is a nice convention for products, since we get
$$
L^p_{\delta_1} \cdot L^q_{\delta_2} \subset L^r_{\delta_1+\delta_2},
\quad \frac{1}{r}=\frac{1}{p}+\frac{1}{q}.
$$

We remark that the Sobolev space above for functions has a simple
interpretation on the conformal half-cylinder with metric
$$ \frac{|dz|^2}{|z|^2} = \frac{dr^2}{r^2}+d\theta^2 = dt^2+d\theta^2 $$
with $t=-\ln r$. Indeed, for a function $f$, the condition
$f\in\hat{L}^{k,p}_\delta$ is equivalent to
$$ e^{\delta t}f \in L^{k,p}\bigg(\frac{|dz|^2}{|z|^2}\bigg), $$
which is the standard weighted Sobolev space on the cylinder.

From this interpretation one easily deduces the following facts
for a function $f$ on the disk, see for example \cite[section 1]{Bi91}:
\begin{enumerate}
\item[(i)] \emph{Sobolev embedding:} for $\frac{1}{2}\geq\frac{1}{p}-\frac{1}{r}$
(with strict inequality for $r=\infty$), one has
\begin{equation}\label{sob_emb}
\hat{L}^{1,p}_{\delta-1}\hookrightarrow L^r_\delta;
\end{equation}
in particular, for $p>2$, one has
$\hat{L}^{1,p}_{\delta-1}\hookrightarrow C^0_\delta$;
\item[(ii)] \emph{control of the function by its radial derivative:}
if $\delta<0$ and $f$ vanishes on the boundary, or if $\delta>0$
and $f$ vanishes near the origin, then
\begin{equation}\label{df_fr}
\|\frac{\partial f}{\partial r}\|_{L^p_{\delta-1}} \geq c \|f\|_{L^p_{\delta}};
\end{equation}
in particular the estimate (\ref{df_fr}) means that
\begin{equation}\label{lhl}
\hat{L}^{1,p}_{\delta-1}=L^{1,p}_{\delta-1} \quad \textrm{ if }\delta<0 ;
\end{equation}
\item[(iii)] for $p>2$ and $\delta>0$, the condition $df\in L^p_{\delta-1}$
implies $f\in C^0$ and then, by applying (\ref{df_fr}):
$$
\|f-f(0)\|_{\hat{L}^{1,p}_{\delta-1}} \leq c \|df\|_{L^p_{\delta-1}}.
$$
\end{enumerate}

From (\ref{lhl}), connection forms
$a\in L^{1,2}_{-2+\delta}(\Omega^1\otimes(\End(E)_0\oplus \End(E)_1))$
actually belong to
$\hat{L}^{1,2}_{-2+\delta} \subset L^p_{-1+\delta}$
for any $p$ by (\ref{sob_emb}).
Now for a gauge transformation $g$, again restricting to the component 
$\End(E)_0\oplus\End(E)_1$, the condition
$Dg g^{-1}\in L^{1,2}_{-2+\delta}\subset L^p_{-1+\delta}$
implies by (iii) that $g$ is continuous, with $g(0)\in
\End(E)_0$, and $g-g(0)\in C^0_\delta\cap L^p_\delta$.

Similarly, an infinitesimal gauge transformation $u$ in
$L^{2,2}_{-2+\delta}(\End(E)_0\oplus\End(E)_1)$ has a well-defined
value $u(0)\in\End(E)_0$, and 
\begin{equation}\label{hL22}
u-u(0)\in\hat{L}^{2,2}_{-2+\delta}.
\end{equation}

We have to generalize this picture from the regular singularity case
to the irregular singularity case. The new ingredient is that now for
sections of $\End(E)_k$ the weight $1/r^k$ is no longer ``equivalent''
to a derivative as in (\ref{df_fr}), and this implies that the
behavior of the weights in the Sobolev embedding $L^{1,2}\hookrightarrow
L^p$ is more involved than in (\ref{sob_emb}), because we somehow have
to separate what comes from the bound on the derivative from what
comes from the bound on the tensor itself.

\begin{lemma}\label{l_sob1}
For $p>2$, one has the Sobolev injections
\begin{eqnarray*}
L^{1,2}_\delta(\End(E)_k) &\hookrightarrow& 
L^p_{\delta+2k/p+(1-2/p)}(\End(E)_k), \\
L^{1,p}_\delta(\End(E)_k) &\hookrightarrow& 
C^0_{\delta+k-2(k-1)/p}(\End(E)_k).
\end{eqnarray*}
\end{lemma}
\begin{proof}
Take $f$ in $L^{1,2}_\delta(\End(E)_k)$, then
$$
\|f\|^2_{L^{1,2}_\delta} =
\int \left( \left| \frac{f}{r^{1+\delta+k}} \right|^2
          + \left| \frac{\nabla f}{r^{1+\delta}} \right|^2
     \right) |dz|^2.
$$
Because of Kato's inequality $|\nabla f|\geq|d|f||$, we can restrict
to the case where $f$ is function, and replace $\nabla f$ by $df$.
With respect to the metric $\frac{|dz|^2}{|z|^{2k}}$, the
above norm transforms into
$$
\|f\|^2_{L^{1,2}_\delta} =
\int \left( \left| \frac{f}{r^{1+\delta}} \right|^2
          + \left| \frac{df}{r^{1+\delta}} \right|^2
     \right) \frac{|dz|^2}{|z|^{2k}} ,
$$
which is equivalent to the $L^{1,2}$-norm
$$
\int (|g|^2+|dg|^2) \frac{|dz|^2}{|z|^{2k}}
$$
of  $g=r^{-1-\delta}f$.
The metric $\frac{|dz|^2}{|z|^{2k}}$ is flat, actually
$$
 \frac{|dz|^2}{|z|^{2k}} = |du|^2, \quad u=\frac{-1}{(k-1)z^{k-1}}.
$$
The problem here is that $z\rightarrow u$ is a $(k-1)$-covering
$\Delta-\{0\} \rightarrow \bC-\bar{\Delta}$: this means that $f$ must
be interpreted on $\bC$ as a section of a rank $(k-1)$ flat unitary bundle.
Nevertheless, still using Kato's inequality, we can apply the standard
Sobolev embedding on $\bC$ to deduce that
$$
\bigg( \int (|g|^2+|dg|^2) |du|^2 \bigg)^\frac{1}{2}
\geq c \bigg( \int |g|^p |du|^2 \bigg)^\frac{1}{p},
$$
that is
$$
\|f\|_{L^{1,2}_\delta} \geq c \|r^{1+\delta-\frac{2k}{p}}f\|_{L^p_{|dz|^2}}
$$
which is exactly the first statement of the lemma.
The proof for the second statement is similar.
\end{proof}

\begin{remark}\label{r_cpct}
Since on a compact manifold these Sobolev embeddings are compact,
it is easy to deduce that for $\delta'<\delta$, the Sobolev embeddings
\begin{eqnarray*}
L^{1,2}_\delta(\End(E)_k) &\hookrightarrow& 
L^p_{\delta'+2k/p+(1-2/p)}(\End(E)_k), \\
L^{1,p}_\delta(\End(E)_k) &\hookrightarrow& 
C^0_{\delta'+k-2(k-1)/p}(\End(E)_k),
\end{eqnarray*}
are compact.
\end{remark}

\begin{remark}
The covering $z^{-(k-1)}=u$ can be used to extend to $L^p$-spaces the
$L^2$-theory which will be done in section \ref{s_analysis}.
\end{remark}

\begin{corollary}\label{c_sob2}
For $k>0$ and $\delta+k+1>0$, one has the injection
$$
L^{2,2}_\delta(\End(E)_k) \hookrightarrow C^0_{\delta'+k+1}
$$
for any $\delta'<\delta$.
\end{corollary}
\begin{proof}
By the previous lemma, if $f\in L^{2,2}_\delta(\End(E)_k)$, then
$$ \nabla f\in L^{1,2}_\delta(\End(E)_k) \subset L^p_{\delta+k-(k-1)(1-2/p)}. $$
Take $p>2$ close enough to $2$ so that $\delta'=\delta-(k-1)(1-2/p)$,
then we get $\nabla f\in L^p_{\delta'+k}$. If $k+1+\delta'>0$, this implies
$$ f=f(0)+f', \quad f'\in C^0_{\delta'+k+1}.$$
Because $f\in L^{2,2}_\delta(\End(E)_k)$, we see that $f(0)=0$.
\end{proof}

We now have the tools to prove that the spaces defined in section
\ref{s_deformations} give us a nice gauge theory.
\begin{lemma}\label{l_alg}
The spaces $L^{2,2}_{-2+\delta}(\End E)$ and $L^{1,p}_{-1+\delta}(\End E)$
(for $p>2$) are algebras, and $L^{1,2}_{-2+\delta}(\End E)$ is a module
over both algebras. 
\end{lemma}
\begin{proof}
We first prove that $L^{1,p}_{-1+\delta}(\End E)$ is an algebra.

We use the fact that $F_k=\oplus_{i\leq k}\End(E)_i$ is a filtration
of $\End E$ by algebras.
Take $u\in L^{1,p}_{-1+\delta}(\End(E)_k)$ and $v\in
L^{1,p}_{-1+\delta}(F_{k-1})$, we want to prove that
$w=uv\in L^{1,p}_{-1+\delta}$; as $F_k$ is an algebra, it is
sufficient to prove that is $w/r^k\in L^p_{\delta-1}$ and
$Dw\in L^p_{\delta-1}$; the statement on $w/r^k$ is clear, since $v\in
C^0$ and $u/r^k\in L^p_{\delta-1}$; now $Dw = (Du)v+u(Dv)$,
but $u,v\in C^0$ and $Du,Dv\in L^p_{-1+\delta}$ implies $Dw\in L^p_{\delta-1}$.

The other statements are proven in a similar way.
\end{proof}

\begin{proof}[Proof of lemma \ref{l_gt}]
The curvature of a connection $D+a\in\cA$ is $F=Da+a\wedge a$. From
the definition of $\cA$, it is clear that $Da\in L^2_{-2+\delta}$; on
the other hand, $a\in L^{1,2}_{-2+\delta}\subset L^4_{-1+\delta}$,
therefore $a\wedge a\in L^2_{-2+2\delta}\subset L^2_{-2+\delta}$, so
$F\in L^2_{-2+\delta}$.

We want to analyze the condition $Dg g^{-1}\in L^{1,2}_{-2+\delta}$
defining the gauge group. Recall that we have a decomposition
$D=D^++\phi$, and the fact that $g$ is unitary actually implies
that both $D^+g g^{-1}$ and $[\phi,g]g^{-1}$ are in
$L^{1,2}_{-2+\delta}\subset L^p_{-1+\delta}$ for any $p>2$.
The condition on $[\phi,g]$ implies $g_k\in L^p_{k-1+\delta}$, and
therefore the condition on $Dg$ implies $dg\in L^p_{-1+\delta}$;
finally, we get $g\in L^{1,p}_{-1+\delta}(\End E)$, and by
lemma \ref{l_sob1}, $$g_k\in C^0_{\delta+k-1-2(k-1)/p}(\End(E)_k).$$
By lemma \ref{l_alg}, we finally deduce that $D^+g$ and $[\phi,g]$ are
in $L^{1,2}_{-2+\delta}(\Omega^1\otimes\End E)$, which implies
$g\in L^{2,2}_{-2+\delta}(\End E)$.

It is now clear that $\cG$ is a Lie group with Lie algebra
$L^{2,2}_{-2+\delta}(\fu(E))$. From lemma \ref{l_alg}, it is easy to
prove the other statements in lemma \ref{l_gt}.
\end{proof}

\section{Analysis on the disk}\label{s_analysis}

In this section, we give some tools to handle the analysis of our
operators with strongly singular coefficients. In order to remain as
elementary as possible, we restrict to $L^2$-spaces, which is
sufficient for our purposes. Moreover, we have
to be careful about the dependence of the constants with respect to
homotheties of the disk, since this is crucial for the compactness
result that we will need later. 

\begin{lemma}\label{l_coercive}
For a $p$-form $u$ with values in $\End E$, with compact support in
$\Delta-\{0\}$, one has
$$
\int |Du|^2+|D^*u|^2 = \int |\nabla u|^2+|\phi\otimes u|^2.
$$
\end{lemma}
\begin{proof}
Integrate by parts the formula \cite[theorem 5.4]{Bi97}
$$
D^*D+DD^*=\nabla^*\nabla+(\phi\otimes)^*\phi\otimes .
$$
\end{proof}

In particular, we apply this formula to get the following consequence.
\begin{corollary}\label{c_control}
If we have a 1-form $u\in L^{1,2}_{-2+\delta}$,
vanishing on $\partial\Delta$, then
$$
\|Du\|_{L^2_{-2+\delta}}+\|D^*u\|_{L^2_{-2+\delta}} \geq
c \big( \|\nabla u\|_{L^2_{-2+\delta}}
      + \|\phi\otimes u\|_{L^2_{-2+\delta}} \big).
$$
On the $\End(E)_k$ part for $k\geq 2$, the estimate is valid for all
weights $\delta$.
For $k=0$ or $1$, the estimate holds only for $\delta>0$ sufficiently small.

Finally, the same estimate holds with the \emph{same} constant $c$ if
we replace $D$ by $h_\varpi^*D$, where $h_\varpi$ is the homothety
$h_\varpi(x)=\varpi x$ in $\Delta$, for $\varpi<1$.
\end{corollary}
\begin{proof}
First consider $u$ section of $\End(E)_k$ for $k\geq 2$.
For a positive function $\rho(r)$ to be fixed later, one has
by lemma \ref{l_coercive},
$$
\|(D+D^*)(\rho^{1-\delta}u)\|^2_{L^2} =
\|\nabla(\rho^{1-\delta}u)\|^2_{L^2} 
+ \|\rho^{1-\delta}\phi\otimes u\|^2_{L^2}.
$$
On the other hand,
$$
\big|[D+D^*,\rho^{1-\delta}]u \big| + \big| [\nabla,\rho^{1-\delta}]u \big| 
\leq c |\rho'\rho^{-\delta}u|.
$$
From these two estimates, we deduce
\begin{eqnarray*}
\|\rho^{1-\delta}(D+D^*)u\|_{L^2}
&\geq& \|(D+D^*)(\rho^{1-\delta}u)\|_{L^2}
      -\|[D+D^*,\rho^{1-\delta}]u\|_{L^2} \\
&\geq& c \big( \|\rho^{1-\delta}\nabla u\|_{L^2}
               +\|\rho^{1-\delta}\phi\otimes u\|_{L^2} \big)
       - c' \|\rho'\rho^{-\delta}u\|_{L^2} .
\end{eqnarray*}
Since $k\geq 2$, we have an (algebraic) estimate
$$ |\phi\otimes u|\geq \lambda_k \frac{|u|}{r^k}, $$
where $\lambda_k$ is the smallest modulus of nonzero eigenvalues of $A_k$.
Now choose $\rho$ in the following way :
\begin{eqnarray*}
\rho(r)=r &&\textrm{for } r\leq \frac{\varepsilon}{2},\\
\rho(r)=\frac{3}{4}\varepsilon &&\textrm{for }r\geq\varepsilon,\\
0\leq\rho'\leq 1.&&
\end{eqnarray*}
Using the fact that $\frac{\rho}{r^2}\geq \frac{1}{2\varepsilon}$ for
$r\leq\varepsilon$, we get the estimate
$$
\|\rho'\rho^{-\delta}u\|^2_{L^2} 
\leq 2\varepsilon \big\|\rho^{1-\delta}\frac{u}{r^2}\big\|^2_{L^2} 
\leq \frac{2\varepsilon}{\lambda_k} \|\rho^{1-\delta}\phi\otimes u\|^2_{L^2}
$$
hence
$$ \|\rho^{1-\delta}(D+D^*)u\|_{L^2} \geq 
\big(c-\frac{2\varepsilon c'}{\lambda_k}\big)
\left( \|\rho^{1-\delta}(D+D^*)u\|_{L^2}
      +\|\rho^{1-\delta}\phi\otimes u\|^2_{L^2} \right). $$
Choosing $\varepsilon$ small enough we get finally 
$$ \|\rho^{1-\delta}(D+D^*)u\|_{L^2} \geq \frac{c}{2}
\left( \|\rho^{1-\delta}(D+D^*)u\|_{L^2}
      +\|\rho^{1-\delta}\phi\otimes u\|^2_{L^2} \right). $$
As $\rho$ coincides with $r$ near $0$, the norm 
$\|\rho^{1-\delta}\cdot\|_{L^2}$ is equivalent to the $L^2_{-2+\delta}$
norm, and the corollary is proven. If $D$ is transformed into
$h_\varpi^*D$, then $\lambda_k$ becomes $\varpi^{1-k}\lambda_k$ which
is bigger, so the estimate still holds ($c$ and $c'$ do not depend on $D$).

In the cases $k=0$ or $1$, one can prove the estimate directly, but
this is a bit more complicated. Another way to prove the lemma is to observe
that in this case, by a conformal change, the operator $D+D^*$ is
transformed into a constant coefficient operator on the
cylinder with metric $\frac{|dz|^2}{|z|^2}$,
and it is clear that it has no kernel. It then follows that
if $\delta$ is not a critical weight (see remark \ref{r_LM}),
the existence of the estimate follows from general
elliptic theory for operators on the cylinder.
Also, the homothety $h_\varpi$ leaves invariant the operator for $k=0$
or $1$, and this implies that the constant does not change.
\end{proof}

\begin{corollary}\label{c_decay}
For $k\geq 2$, if we have a $L^2_\delta$-solution $u$ with values in
$\End(E)_k$ of the equation $(D+D^*)u=0$, then $|u|=O(r^\gamma)$ for
any real $\gamma$. 
\end{corollary}
\begin{proof}
Let $\chi_\epsilon(r)$ be a cut-off function, so that $\chi_\epsilon u=u$
for $\epsilon<r<1/2$. We can choose $\chi_\epsilon$ so that
$$|d\chi_\epsilon|\leq \frac{c}{r}.$$
Then
$$ |(D+D^*)(\chi_\epsilon u)|\leq \frac{c}{r} |u|. $$
By corollary \ref{c_control} this gives
\begin{eqnarray*}
\|\nabla(\chi_\epsilon u)\|_{L^2_{\delta-1}}
+\|\chi_\epsilon u\|_{L^2_{\delta+k-1}}
&\leq& c'\|(D+D^*)(\chi_\epsilon u)\|_{L^2_{\delta-1}} \\
&\leq& cc'\|u\|_{L^2_\delta}.
\end{eqnarray*}
Taking $\varepsilon\rightarrow 0$, we get $u\in L^2_{\delta+k-1}$ and
$\nabla u\in L^2_{\delta-1}$.
Since $k>1$, we can iterate the argument and deduce that $u$ and
$\nabla u$ (hence $d|u|$) belong to $L^2_\gamma$ for any $\gamma$.
The corollary follows.
\end{proof}

\begin{lemma}\label{l_disk}
On the disk, the Laplacian
$$D^*D+DD^*:L^{2,2}_{-2+\delta}(\Omega^i\otimes\End E)\rightarrow
            L^2_{-2+\delta}(\Omega^i\otimes\End E)$$
with Dirichlet condition on the boundary, is an isomorphism for small
weights $\delta>0$.

If we restrict to the components $\End(E)_n\oplus\cdots\oplus\End(E)_2$,
then the same holds for any weight $\delta$.
\end{lemma}
\begin{proof}
Begin by the components in $\End(E)_k$ for $k\geq 2$. We claim that a
solution of $(DD^*+D^*D)u=v$ is obtained by minimizing the functional
$$
\int \frac{1}{2}\big(|Du|^2+|D^*u|^2\big) - \langle u,v \rangle
$$
among $u\in L^{1,2}_{-1}(\End(E)_k)$ vanishing on the boundary; indeed,
by lemma \ref{l_coercive}, we have for such $u$
\begin{equation}\label{332}
\int \big(|Du|^2+|D^*u|^2\big) \geq c \int \frac{|u|^2}{r^{2k}},
\end{equation}
so the functional is coercive and a solution can be found; note that
it is sufficient to have $r^k v\in L^2$. The solution satisfies
an equation
$$
(dd^*+d^*d)u=v+P_0\bigg(\frac{u}{r^{2k}}\bigg)+
               P_1\bigg(\frac{\nabla u}{r^k}\bigg),
$$
where $P_0$ and $P_1$ are bounded algebraic operators,
so by elliptic regularity $u\in L^{2,2}_{-k-1}$. In conclusion
we get, for the Dirichlet boundary condition, an isomorphism
\begin{equation}\label{334}
DD^*+D^*D:L^{2,2}_{-k-1}\rightarrow L^2_{-k-1}.
\end{equation}

Now we want to prove that this $L^2$-isomorphism actually extends to all
weights: we proceed by proving that the $L^2$-inverse is
continuous in the other weighted spaces. 
For any weight $\gamma$, it is sufficient to prove an estimate
$$
\|\rho^\gamma(DD^*+D^*D)u\|_{L^2_{-k-1}} \geq c
\|\rho^\gamma u\|_{L^{2,2}_{-k-1}},
$$
where $\rho$ is some function which coincides with $r$ near $0$, as in
the proof of corollary \ref{c_control}. As in the proof of this corollary,
the estimate is deduced from a control on the commutator
$[DD^*+D^*D,\rho^\gamma]$, obtained after a careful choice of $\rho$. 
The details are left to the reader.

Now let us look at the component $\End(E)_1$. Actually the estimate
(\ref{332}) still holds, therefore the isomorphism (\ref{334}) is
true, and remains true for small perturbations $-2+\delta$ of the
weight $-2$.

Finally, for the component $\End(E)_0$, we simply have the usual
Laplacian $dd^*+d^*d$ on the disk to study in weighted Sobolev spaces,
so it is a classical picture: the inverse is the usual
solution of the Dirichlet problem on the disk.
The question is to check the regularity in the weighted Sobolev spaces.
It is useful to proceed in the following way, using general theory for
elliptic operators on cylinders: 
the operator 
\begin{equation}\label{ddL}
dd^*+d^*d:\hat{L}^{2,2}_{-2+\delta}\rightarrow L^2_{-2+\delta}
\end{equation}
translates on the cylinder into the operator 
\begin{equation}\label{ddLC}
-\frac{\partial^2}{\partial t^2}-\frac{\partial^2}{\partial\theta^2}:
e^{\delta t}L^{2,2}\longrightarrow e^{\delta t}L^2 .
\end{equation}
Note that the operator for $k=1$ is just the same, with an additional
term $\lambda_1^2$, where $\lambda_1$ is the eigenvalue of the action
of $A_1$.

Here the weight $\delta=0$ is critical, because of the solutions
$a+bt$ in the kernel of (\ref{ddLC}). Nevertheless, for small
$\delta\neq 0$ the operator (\ref{ddLC}) becomes Fredholm, and because
the operator is selfadjoint, we get by formula \cite[theorem
7.4]{LoMO85} the index $+1$ for small negative $\delta$ and $-1$ for
small positive $\delta$. Coming back on the disk, this means the 
operator (\ref{ddL}) has index $-1$. As
$\hat{L}^{2,2}_{-2+\delta}\subset L^{2,2}_{-2+\delta}$ has codimension
1 by (\ref{hL22}), this means that
$dd^*+d^*d:L^{2,2}_{-2+\delta}\rightarrow L^2_{-2+\delta}$ 
has index 0. Therefore it is an isomorphism.
\end{proof}

\begin{remark}\label{r_LM}
On the $\End(E)_0\oplus\End(E)_1$ part, the theory of elliptic
operators on cylinders gives also information for all weights $\delta$.
Namely, the problem is Fredholm if $\delta$ avoids a discrete set of
critical weights (corresponding to the existence of solutions
$(a+bt)e^{-\delta t}$).
This easily follows from the following fact: 
if one considers the problem (\ref{ddLC}) on the whole cylinder
(or, equivalently, $DD^*+D^*D$ on $\bR^2-\{0\}$), then it is an
isomorphism outside these critical weights.
\end{remark}

Finally, we deduce the decay of the solutions of the selfduality equations.

\begin{lemma}\label{l_sol_decay}
If we have on the disk a solution 
$a\in L^{1,2}_{-2+\delta}(\Omega^1\otimes\End E)$ of an
equation $(D+D^*)a=a\odot a$, then in the decomposition $a=\sum a_k$
we have the following decay for $a$:
\begin{enumerate}
\item if $k\geq 2$, then $|a_k|=O(r^\gamma)$ for any $\gamma$;
\item if $k=0$ or $1$, then $|a_k|=O(r^{-1+\delta})$.
\end{enumerate}
\end{lemma}
\begin{proof}
We have $L^{1,2}_{-2+\delta}\subset L^p_{-1+\delta}$ for any $p$. On the
other hand, for $q>2$ close enough to $2$, one has the inclusion 
$L^{1,2}_{-2+\delta}(\End(E)_k)\subset L^q_{-2+k}$. 
Now take $p$ big enough so that $1/p+1/q=1$, we get
$a\odot a_k\in L^2_{-3+k+\delta}$.
In particular, we deduce that $(a\odot a)_k\in L^2_{-3+k+\delta}$ for all
$k\geq 2$. It is now easy to adapt the proof of the corollary \ref{c_decay} to
get $a_k\in L^2_{-3+2k+\delta}$ for all $k\geq 2$, and therefore
$a_k\in L^{1,2}_{-3+\delta+k}\subset L^{1,2}_{-1+\delta}$.

Iterating this, we get that for $k\geq 2$ one has $a_k\in L^{1,2}_\gamma$
for any $\gamma$, and we deduce the bound on $a_k$.

For the $b=a_0+a_1$ part, we can write the problem as
$$
(D+D^*)b=b\odot b + \textrm{small perturbation},
$$
with an initial bound $b\in L^{1,2}_{-2+\delta}$, and therefore $b\odot
b\in L^p_{-2+2\delta}$ for any $p>2$.
This is now a problem which translates into a constant coefficient
elliptic problem on the conformal cylinder, so that elliptic regularity
gives at once  $b\in L^{1,p}_{-2+\delta}\subset C^0_{-1+\delta}$.
\end{proof}

\section{Moduli spaces}\label{s_moduli}
Consider now a compact Riemann surface $X$ with finitely many
marked points $p_i$, and a complex vector bundle $E$ over $X$, with a
hermitian metric $h$. Choose an initial connection $D_0$ on $E$, such
that in some unitary trivialization of $E$ around each $p_i$, the
connection $D_0$ coincides with the local model (\ref{f_modDR}).
Of course on the interior of $X-\{p_i\}$, the connection $D_0$ is not
flat in general.

Define $r$ to be a positive function which coincides with $|z|$ around
each puncture, then we can define global Sobolev spaces on $X$ as in
section \ref{s_deformations}, and therefore a space of connections
$\cA=D_0+L^{1,2}_{-2+\delta}(\Omega^1\otimes\End E)$, and a gauge group
$\cG$ as in (\ref{d_gauge}) acting on $\cA$. 

The result \ref{l_disk} on the disk now implies globally.
\begin{lemma}\label{l_fredholm}
If $A\in\cA$, then the operator
$$D_A^* D_A:L^{2,2}_{-2+\delta}(\fu(E)) \rightarrow L^2_{-2+\delta}(\fu(E))$$
is Fredholm, of index 0.
\end{lemma}
\begin{proof}
First it is sufficient to prove that the Laplacian $D_0^*D_0$ is
Fredholm, because the 1-form $a\in L^{1,2}_{-2+\delta}(\Omega^1\otimes\End E)$
gives only a compact perturbation, see remark \ref{r_cpct}.
For $D_0^*D_0$, we can glue the inverse coming from lemma \ref{l_disk}
near the punctures with a parametrix in the interior, this gives a
parametrix which is an exact inverse near infinity, implying that the
operator is Fredholm.

Let $i_\delta$ be the index of the slightly different operator,
$$ P_\delta=D_0^*D_0:\hat{L}^{2,2}_{-2+\delta}(\fu(E))
 \rightarrow L^2_{-2+\delta}(\fu(E)) ;$$
this means that we now do not allow nonzero values at the origin for
the $\End(E)_0$ part. We claim that
\begin{eqnarray}
& & i_\delta=-i_{-\delta}, \label{ind_opp} \\
& & i_\delta-i_{-\delta}=-2\dim\fu(E)_0. \label{ind_dif}
\end{eqnarray}
From these two assertions, it follows immediately that
$$ i_{\delta}=-\dim\fu(E)_0,$$
and therefore the index of the initial operator is $0$.

Now let us prove first (\ref{ind_dif}). We have to calculate the
difference between the indices of $P_\delta$ and $P_{-\delta}$.
The operator does not change in the interior of $X$, so by the
excision principle, the difference comes only from what happens
at the punctures, and it is sufficient to calculate it for the model
Dirichlet problem: this has been done in lemma \ref{l_disk} and its proof.

Now we prove (\ref{ind_opp}). This comes from formal $L^2$-selfadjointness
of $P_\delta$: observe that the dual of $L^2_{-2+\delta}$ is
identified to $L^2_{-\delta}$, therefore the cokernel of $P_\delta$
consists of solutions $u\in L^2_{-\delta}$ of the equation $D_0^*D_0u=0$.
The behavior of such a $u$ comes from lemma \ref{l_disk}: near the
punctures, the components $u_k$ for $k\geq 2$ decay quicker than
any $r^\gamma$. This fact combined with elliptic regularity implies
$u\in\hat{L}^{2,2}_{-2-\delta}(\fu(E))$, that is $u\in\ker
P_{-\delta}$. Therefore $\coker P_\delta=\ker P_{-\delta}$.
\end{proof}

We want to consider the quotient space $\cA/\cG$. If $A\in\cA$ is
irreducible, then $D_A$ has no kernel and the cokernel of $D_A$ is
simply the kernel of $D_A^*$. From this and the lemma, it is classical
to deduce that the irreducible part $\cA^{irr}/\cG$ of the quotient is
a manifold, with tangent space at a connection $A$ given by
\begin{equation}
T_{[A]}(\cA^{irr}/\cG)=
\{a\in L^{1,2}_{-2+\delta}(\Omega^1\otimes\End E), \im(D_A^*a)=0\}.
\end{equation}

The moduli space $\cM\subset\cA/\cG$ we consider is defined by the equations
\begin{equation}\label{f_FG}
F_A=0, \quad G_A=0 .
\end{equation}
These equations are not independent, since writing $A = A^+ + \phi_A$,
we have the decomposition of $F_A$ and $G_A$ into selfadjoint and
antiselfadjoint parts given by
\begin{eqnarray*}
F_A &=& D_{A^+}\phi_A + (F_A+\phi_A\wedge \phi_A),\\
\Lambda G_A &=& -\Lambda D_{A^+}\phi_A + iD_{A^+}^*\phi_A.
\end{eqnarray*}
Therefore, as is well-known, the equations (\ref{f_FG}) are equivalent to
\begin{equation}\label{f_FG2}
F_A=0, \quad D_{A^+}^*\phi_A=0 .
\end{equation}
The space $\cA$ is a flat hyperK{\"a}hler space, for the
standard $L^2$-metric, and with complex structures $I$, $J$ and $K=IJ$
acting on $a\in\Omega^1\otimes\End(E)$ by
\begin{equation}\label{f_IJ}
I(a)=ia,\quad J(a)=i(a^{0,1})^*-i(a^{1,0})^*.
\end{equation}
The complex structure $I$ is the natural complex structure on
connections, and $J$ is the natural complex structure on Higgs bundles.

Hitchin \cite{Hi87} observed that
the equations (\ref{f_FG2}) are the zero set of the hyperK{\"a}hler
moment map of the action of $\cG$ on $\cA$.
The linearization of the equations (\ref{f_FG2}) is simply
\begin{equation}\label{f_linFG}
D_Aa=0,\quad \re(D_A^*a)=0.
\end{equation}
If $[A]\in\cM$, there is an elliptic deformation complex governing the
deformations of $A$:
$$
L^{2,2}_{-2+\delta}(\fu(E)) \stackrel{D_A}{\rightarrow}
L^{1,2}_{-2+\delta}(\Omega^1\otimes\End E) \stackrel{D_A+D_A^*}{\longrightarrow}
L^2_{-2+\delta}((\Omega^2\otimes\End E)\oplus i\fu(E))
$$

As in the proof of lemma \ref{l_fredholm}, the Laplacians of the
complex are Fredholm, with index 0, and we get the following result.
\begin{lemma}\label{l_complex}
The cohomology groups of the deformation complex are finite
dimensional.\qed
\end{lemma}

\begin{lemma}\label{l_beh}
If $A\in\cM$, then there is a gauge in which near a puncture,
$A=D_0+a$, and $a$ decays as in the conclusion of lemma \ref{l_sol_decay}.
\end{lemma}
\begin{proof}
We can write globally $A=D_0+a$. Let $\chi_\epsilon$ be a cutoff
function, such that
\begin{enumerate}
\item $\chi_\epsilon=1$ in a disk of radius $\epsilon$ near each puncture;
\item $\chi_\epsilon=0$ outside the disks of radius $2\epsilon$ near
each puncture
\item $|d\chi_\epsilon|\leq c/r$ for some constant $c$.
\end{enumerate}
Consider the connections $A_\epsilon=D_0+\chi_\epsilon a$. 
This is a continuous path of connections in $\cA$, converging to $D_0$.
\begin{claim}[Coulomb gauge]
For $\epsilon$ sufficiently small, there exists a gauge transformation
$g_\epsilon\in \cG$, such that
\begin{equation}\label{coulomb}
\im\big(D_0^*(g_\epsilon(A_\epsilon)-D_0)\big)=0.
\end{equation}
\end{claim}
Suppose the claim is proven, fix some $\epsilon$ for which we have a
Coulomb gauge; since $A=A_\epsilon$ in a disk of radius $\epsilon$
near the punctures, this means that, restricting to this disk,
$$
g_\epsilon(A)=D_0+a,\quad \im D_0^*a=0.
$$
Using this condition, the equations (\ref{f_FG2}), with linearization
(\ref{f_linFG}),  can be written $(D_0+D_0^*)a=a\odot a$ and the
result follows from lemma \ref{l_sol_decay}.

It remains to prove the claim. We try to find $g_\epsilon=e^u$ solving
the equation (\ref{coulomb}), with $u\in L^{2,2}_{-2+\delta}(\fu(E))$.
The linearization of the problem is
$$ D_0^*D_0u=\im D_0^*(\chi_\epsilon a). $$
But the operator $D_0^*D_0$ is Fredholm of index 0 by lemma
\ref{l_fredholm}. If it is invertible then by the implicit function
theorem, we get the solution $g_\epsilon$ we wanted. If not, we still
get an isomorphism after restricting to the space $\ker(D_0^*D_0)^\perp$;
fortunately, the operator 
$\im D_0^*(e^u(A)-D_0)$ takes its values in the same space, so we can
still apply the implicit function theorem after restricting to it.
\end{proof}

If $A$ is irreducible, then the kernel of $D_A$ on $\End(E)$ vanishes,
and is equal to the kernel of $D_A^*$ on $\Omega^2\otimes\End E$.
From lemma \ref{l_complex}, it now follows that the equations
(\ref{f_FG}) are transverse, and we therefore get the following result.
\begin{theorem}\label{t_m}
The moduli space $\cM^{irr}$ is a smooth hyperK{\"a}hler manifold, with
tangent space at $A$ given by
$$
T_{[A]}\cM=\{a\in L^2(\Omega^1\otimes\End E), D_Aa=0, D_A^*a=0\}
          =L^2H^1(\End E).
$$
The metric is the natural $L^2$-metric.
Under the secondary assumption ($A_1$ is semisimple), the metric is complete
if the moduli space $\cM$ does not contain reducible points.
\end{theorem}
\begin{proof}
First we prove that $L^2$-cohomology calculates the $H^1$ of the
elliptic complex.
We have to prove that a $L^2$-harmonic 1-form actually belongs to the
space $L^{1,2}_{-2+\delta}$.
This is the infinitesimal version of lemma \ref{l_sol_decay} and is
even simpler to prove (and because the equations are conformally
invariant, one can do the local calculations on the disk with respect
to the flat metric).

The difficult point is to prove that the metric is complete.
Suppose we have a geodesic curve $([A_t])$ in $\cM$ parameterized by arclength,
of finite length $\ell$. We want to extend it a bit.
We can lift it to a horizontal curve $(A_t=A_0+a_t)$, hence
$\dot{a}_t$ satisfies 
\begin{eqnarray}
 (D_{A_t}+D_{A_t}^*)\dot{a}_t &=& 0, \label{C1} \\ 
 \int_X |\dot{a}_t|^2 &=& 1. \label{C2}
\end{eqnarray}
These two equalities do not depend on the metric on $X$.
If we choose a metric $g$, then, 
decomposing $\nabla_{A_t}=\nabla_t^++\phi_t$,
we get the Weitzenb{\"o}ck formula \cite[theorem 5.4]{Bi97}:
$$ (\nabla_t^+)^*\nabla_t^+\dot{a}_t
  + \phi_t^*\phi_t\dot{a}_t
  + \frac{\scal^g}{2} \dot{a}_t = 0. $$
We cannot integrate this equation against $\dot{a}_t$ for a smooth
$g$, because the integral is divergent. Nevertheless, we use the
freedom of the metric to choose $g$ which coincides near each puncture
with the flat metric
$$ \frac{|dz|^2}{|z|^{2(1-\delta')}} $$
for some local coordinate $z$ and some positive $\delta'<\delta$.
Because $\dot{a}_t\in L^{1,2}_{-2+\delta}$ now one can integrate by
parts and get
$$
\int_X \big(
  |\nabla_t^+\dot{a}_t|^2 + |\phi_t\dot{a}_t|^2
+ \frac{\scal^g}{2} |\dot{a}_t|^2
     \big) vol^g = 0.
$$
Here all norms are taken with respect to $g$.
Because the scalar curvature of $g$ is bounded, and using Kato's
inequality, we deduce
\begin{equation}\label{C5}
\int_X |d|\dot{a}_t||^2 vol^g \leq c \int_X |\dot{a}_t|^2 \leq c.
\end{equation}
Since the $L^2$-norm of 1-forms in conformally invariant, we can write
on each disk near a puncture this equality with respect to the flat
metric $|dz|^2$:
\begin{equation}\label{C3}
\int_{r<1} |d|\dot{a}_t|_g|^2 |dz|^2 =
\int_{r<1} | d|r^{1-\delta'}\dot{a}_t| |^2 |dz|^2 \leq c.
\end{equation}
\begin{claim}One has the estimate
$$ \int_{r<1} \left|r^{-1+\varepsilon}f\right|^2 |dz|^2 \leq c
   \bigg( \int_{r<1} |df|^2 |dz|^2 + \int_{\frac{1}{2}<r<1} |f|^2 |dz|^2 
   \bigg). $$
\end{claim}
Using the claim, we deduce from (\ref{C2}) and (\ref{C3}), for some positive
$\delta''<\delta'$, the estimate 
$$
\int_{r<1} \bigg( | r^{1-\delta''} d|\dot{a}_t| |^2 +
                 | r^{-\delta''}\dot{a}_t |^2 \bigg) |dz|^2 \leq c.
$$
In particular, the $L^2$-estimate is now slightly better than (\ref{C2}).
Furthermore, the Sobolev embedding (\ref{sob_emb}) implies
$$ \| \dot{a}_t \|_{L^p_{-1+\delta''}} \leq c. $$
This estimate holds on every disk near the singularities, but also in
the interior of $X$, applying the Sobolev embedding to (\ref{C5});
hence it is now a global estimate on $X$.
Since $a_t=\int_0^t \dot{a}_t dt$, we get also the estimate on $X$:
$$ \| a_t \|_{L^p_{-1+\delta''}} \leq c. $$
Now rewrite the equation (\ref{C1}) in the form
$$
(D_{A_0}+D_{A_0}^*)\dot{a}_t = a_t \odot \dot{a}_t.
$$
Choose $\delta''$ so close to $\delta$ that $2\delta''\geq\delta$.
Then, from the multiplication 
$$L^4_{-1+\delta''}\cdot L^4_{-1+\delta''}\subset
L^2_{-2+2\delta''}\subset L^2_{-2+\delta},$$
we deduce that $\dot{a}_t$ remains bounded in $L^{1,2}_{-2+\delta}$,
which means that $a_t$ has a limit $a_\ell$ in $L^{1,2}_{-2+\delta}$ when $t$
goes to $\ell$.
The limiting $A_0+a_\ell$ is again a solution of Hitchin's
equations, so represents a point of $\cM$. Since there is no reducible
solution, it is a smooth point and the geodesic can be extended.

There remains to prove the claim. It is a simple application of (\ref{df_fr}):
$$ \int |df|^2|dz|^2 
\geq \int |r^\varepsilon df|^2|dz|^2
\geq c \int |r^{-1+\varepsilon}f|^2|dz|^2 $$
if $f$ vanishes at the boundary. If not, then use a cut-off function
$\chi$ so that $\chi=1$ for $r\leq\frac{1}{2}$, and apply the above
inequality to $\chi f$:
$$ \int |df|^2|dz|^2 + \int_{\frac{1}{2}<r<1} |f|^2|dz|^2 
\geq \int |d(\chi f)|^2|dz|^2 
\geq \int |r^{-1+\varepsilon}\chi f|^2|dz|^2.$$
The claim follows.
\end{proof}

\section{Complex moduli spaces and harmonic metrics}\label{s_complex}

There are two complex moduli spaces that we would like to consider. We
still have some reference connection $D_0\in\cA$ and recall that we
can decompose $D_0=D_0^+ + \phi_0$.

We have defined the unitary gauge group $\cG$ by the condition
$D_0gg^{-1}\in L^{1,2}_{-2+\delta}$. This condition implies that both
$D_0^+g g^{-1}$ and $g\phi_0 g^{-1}$ are in $L^{1,2}_{-2+\delta}$.

For complex transformations, this is no longer true, and we have to
define directly the complexified gauge group $\cG_{\bC}$ as the space of
transformations $g\in GL(E)$ such that $D_0^+ g g^{-1}$ and
$g\phi_0 g^{-1}$ belong to $L^{1,2}_{-2+\delta}$.
As in lemma \ref{l_gt}, this definition makes $\cG_{\bC}$ into a Lie
group with Lie algebra 
$$ \Lie(\cG_{\bC})=L^{2,2}_{-2+\delta}(\End E). $$
We now complexify the action of $\cG$ on $\cA$ for both the complex
structures $I$ and $J$.

For the first complex structure, the action is simply given by
$$ D \longrightarrow g\circ D\circ g^{-1}=D-Dg g^{-1} $$
and the associated moduli space is the moduli space of flat
connections on $E$, defined by
\begin{equation}
\cM_{DR,an} = \{ A\in\cA, F_A=0 \}/\cG_{\bC}.
\end{equation}
The DR subscript means ``De Rham'' moduli space, as in Simpson's
terminology, and the ``an'' is for ``analytic'', by contrast with
the moduli space of flat connections $\cM_{DR,alg}$ with some fixed
behavior at the punctures that one can define algebraically.

For the second complex structure, we get a different action, namely
the action of $\cG_{\bC}$ on $\cA$ seen as a space of Higgs bundles:
a connection $A=D_0+a\in\cA$ can be identified by (\ref{f_higgs})
with the Higgs bundle $(\db_A,\theta_A)$, writing
$$
\db_A = \db_0 + \frac{a^{0,1}-(a^{1,0})^*}{2},\quad
\theta_A = \theta_0 + \frac{a^{1,0}+(a^{0,1})^*}{2},
$$
and the action of the complexified gauge group is simply
$$ (\db_A,\theta_A) \longrightarrow
(g\circ\db_A\circ g^{-1},g\theta_A g^{-1}). $$
The associated moduli space is
\begin{equation}
\cM_{Dol,an} = \{ A\in\cA, G_A=0 \}/\cG_{\bC}.
\end{equation}
Again the subscript ``Dol'' stands for Dolbeault moduli space, and the
``an'' is used to distinguish with the algebraic moduli space
$\cM_{Dol,alg}$ of Higgs bundles with Higgs field of fixed polar
part at each puncture.

In the two cases there is a notion of (analytic) stability which leads to
spaces $\cM_{DR,an}^s \subset \cM_{DR,an}$ and $\cM_{Dol,an}^s \subset
\cM_{Dol,an}$. In both cases, stability means that for some class of
subbundles, the slope of the subbundle (the analytic degree defined by
the metric, divided by the rank) must be smaller than the slope of
the bundle. Here we will only recall the definition and refer the
reader to \cite[section 6]{Si90} for details.

In the Dolbeault case, the class of subbundles to consider is the
class of holomorphic $L^{1,2}$-subbundles, that is holomorphic
subbundles $\cF$ (outside the punctures), stable under the Higgs field, and
defined by an orthogonal projection $\pi$ such that $\db_A\pi\in L^2$.
The corresponding analytic degree is obtained by integrating the
curvature of the connection induced on the subbundle by the metric,
\begin{eqnarray*}
\deg^{an}\cF&=&\frac{i}{2\pi}\int \tr(\pi F_A)-|D''_A\pi|^2\\
            &=&\frac{i}{2\pi}\int \tr(\pi F_{A^+})-|\db_A\pi|^2.
\end{eqnarray*}

In the DeRham case, there is a similar picture: one has to consider
flat subbundles defined by an orthogonal projection $\pi$ such that
$D_A\pi\in L^2$. 
On a compact manifold, the degree of a flat subbundle is always $0$,
and stability reduces to semisimplicity, that is there is no flat
subbundle; in the noncompact case, the parabolic structure at the
punctures may have a nonzero contribution to the degree of a flat subbundle.

\begin{theorem}\label{t_he}
Suppose $\deg^{an} E=0$. Then the natural restriction maps
$$
\cM_{Dol,an}^s \longleftarrow \cM^{irr} \longrightarrow \cM_{DR,an}^s
$$
are isomorphisms.
\end{theorem}

We will prove the theorem in section \ref{s_proof}.

\section{The Dolbeault moduli space}\label{s_Dol}

In this section, we prove that elements of the analytic Dolbeault moduli space
$\cM_{Dol,an}$ actually correspond to true meromorphic Higgs bundles,
with fixed parabolic structure and fixed polar part of the Higgs field,
on the Riemann surface $X$. This gives a correspondence between
$\cM_{Dol,an}^s$ and $\cM_{Dol,alg}^s$.

In $X-\{p_i\}$, a $L^{1,2}_{loc}$-$\db$-operator has holomorphic sections
which define a structure of holomorphic bundle. The remaining question
is local near the punctures. Fix the local model $D_0$ around a
puncture as in section \ref{s_model}, with underlying Higgs bundle
given by (\ref{f_modDol}):
$$
\db_0=\db-\frac{1}{2}\re(A_1)\frac{d\zb}{\zb}, \quad 
\theta_0= \frac{1}{2}\sum_1^n A_i \frac{dz}{z^i}
         - \frac{\beta}{2} \frac{dz}{z} .
$$
Now consider a Higgs bundle
$(\db_0+a,\theta_0+\vartheta)\in\cA$. The following lemma enables to
take the $\db$-operator to a standard form. 
\begin{lemma}\label{l_db}
There is a complex gauge transformation $g$, defined in a neighborhood
of the origin, such that
\begin{enumerate}
\item $g$ is continuous, and in the decomposition $g=g_0+\cdots+g_n$, one has
$g_i/r^{i-1+\delta'}$ continuous for some $\delta'<\delta$;
\item $g(\db_0+a)=\db_0$.
\end{enumerate}
\end{lemma}

To prove the lemma, we first need the following statement.
\begin{lemma}\label{l_db1}
Take $\delta\in\bR-\bZ$ and $p>2$. Then the problem 
\begin{equation}\label{e_dbpro}
\frac{\partial f}{\partial\zb}=g
\end{equation}
in the unit disk has a solution $f=T_0g$ such that
\begin{equation}\label{e_dbcon}
\|f\|_{C^0_{-1+\delta}} \leq c \|g\|_{L^p_{-2+\delta}}.
\end{equation}
The same is true, if $\delta-\re\lambda\not\in\bZ$, for the problem
$$ \frac{\partial f}{\partial\zb}+\frac{\lambda}{2} \frac{f}{\zb}=g. $$
\end{lemma}
\begin{proof}
First one can restrict to the case $0<\delta<1$. Indeed, if
$\delta=[\delta]+\delta_0$, then the problem (\ref{e_dbpro}) is equivalent to
$$
\frac{\partial (z^{-[\delta]}f)}{\partial\zb}=z^{-[\delta]}g,
$$
and the wanted estimate (\ref{e_dbcon}) becomes equivalent to the estimate
$$
\|z^{-[\delta]}f\|_{C^0_{-1+\delta_0}} 
\leq c \|z^{-[\delta]}g\|_{L^p_{-2+\delta_0}}.
$$
Hence we may suppose $0<\delta<1$.

Now for this range of $\delta$, we claim that the Cauchy kernel
$$
f(z) = \int \frac{g(u)}{z-u}|du|^2
$$
gives the inverse we need. Indeed, by the H{\"o}lder inequality, we get
$$
|f(z)| \leq \|g\|_{L^p_{-2+\delta}} \left( \int
\frac{|du|^2}{|u|^{2-\delta\frac{p}{p-1}}|z-u|^{\frac{p}{p-1}}}
                                    \right)^{\frac{p-1}{p}}
$$
and the result follows from the estimate
$$
\int \frac{|du|^2}{|u|^{2-\delta\frac{p}{p-1}}|z-u|^{\frac{p}{p-1}}}
\leq \frac{c}{|z|^{(1+\delta)\frac{p}{p-1}}}
$$
which is left to the reader.

For the second problem, observe that
$$ \frac{\partial}{\partial\zb}+\frac{\lambda}{2\zb} =
  r^{-\lambda} \circ \frac{\partial}{\partial\zb} \circ r^{\lambda}. $$
If $\delta-\re\lambda\not\in\bZ$, we can take the inverse
$$ f=r^{-\lambda} T_0 r^\lambda g.$$
\end{proof}

\begin{proof}[Proof of lemma \ref{l_db}]
Define spaces of $u\in\End E$ and $a\in\Omega^1\otimes\End E$ by
\begin{eqnarray*}
U_\delta &=& \{u, u_k\in C^0_{k-1+\delta}, u_0\in C^0_\delta\}, \\
A_\delta &=& \{a, a_k\in L^p_{-2+\delta+k}, a_0\in L^p_{-1+\delta}\}.
\end{eqnarray*}
Observe that since $a\in L^{1,2}_{-2+\delta}$, we have $a_k\in
L^p_{\delta+k-2-\epsilon}$ by lemma \ref{l_sob1}, and therefore $a\in
A_{\delta'}$ for some weight $\delta'<\delta$.

The problem $g(\db_0+a)=\db_0$ can be written (with $g=1+u$)
\begin{equation}\label{711}
\db_0 u - u a = a.
\end{equation}
The operator $\db_0$ is of type studied in lemma \ref{l_db1};
hence we get a continuous right inverse 
$T:A_{\delta'}\rightarrow U_{\delta'}$.
We find a solution $u\in U_{\delta'}$ of (\ref{711}) by a fixed point problem,
looking at a solution $u\in U_{\delta'}$ of 
$$ u = T(ua+a). $$
For this, we need $u\rightarrow T(ua+a)$ to be contractible; but
$$ \|T(ua+a)-T(va+a)\|_U \leq   c\|u-v\|_U \|a\|_A, $$
so this is true if $\|a\|_A$ is small enough.

Let $h_\varpi$ be a homothety taking the disk of radius 1 to the
disk of radius $\varpi$, then it is easy to see that
$$
\|h_\varpi^* a\|_A \leq \varpi^\delta \|a\|_A.
$$
On the other hand, the operator $\db_0$ is unchanged by the homothety
$h_\varpi$, so for $\varpi$ small enough, the operator
$u\rightarrow T(ua+a)$ becomes contractible and we can solve the problem.
\end{proof}

From the lemma, we deduce a basis of holomorphic sections for $\db_A$,
as in (\ref{f_sigma}),
$$
\sigma_i=g^{-1} |z|^{\re\mu_i-[\re\mu_i]} e_i.
$$
This basis defines a holomorphic extension of the bundle $(E,\db_A)$ over the
puncture, which is characterized by the fact the sheaf of holomorphic
sections of this bundle is the sheaf of bounded holomorphic sections
outside the puncture. Moreover, the growth of the holomorphic sections
is the same as the model (\ref{f_alpha}), that is
\begin{equation}\label{m_bhvr}
|\sigma_i| \sim |z|^{\alpha_i}
\end{equation}
with $\alpha_i=\re\mu_i-[\re\mu_i]$, and these different orders of
growth define on the extension a parabolic structure, whose weights
are $\alpha_1,\dots,\alpha_r$.

Finally, in the holomorphic basis $(\sigma_i)$, the Higgs field becomes
$$
\theta=g(\theta_0+\vartheta)g^{-1}=\theta_0+\vartheta',\quad
\vartheta'=g[\theta_0,g^{-1}]+g\vartheta g^{-1},
$$
and we now have simply $\db\vartheta'=0$. From the bounds on $g$ in
the lemma, we deduce that actually $\vartheta'$ is holomorphic,
therefore the polar part of the Higgs field is exactly, as in (\ref{f_theta}),
$$
\theta_0= \frac{1}{2}\sum_1^n A_i \frac{dz}{z^i}
        - \frac{\beta}{2} \frac{dz}{z}.
$$

The above construction means that we have constructed a map
\begin{equation}\label{a_ext}
\cM_{Dol,an} \rightarrow \cM_{Dol,alg}
\end{equation}
by defining a canonical extension. Actually, writing (\ref{a_ext}) is
not completely correct, and it would be better to say that we have a
functor between the two corresponding categories (indeed one can prove
that morphisms in the space $L^{2,2}_{-2+\delta}$ extend to
holomorphic morphisms of the extensions).
This functor is actually an equivalence of categories, because the
arrow (\ref{a_ext}) can be inverted through the following lemma.

\begin{lemma}\label{l_initial}
Let $(\cE,\theta)$ be a meromorphic Higgs bundle on $X$, with
parabolic structure at the punctures having weights $\alpha_1,\dots,\alpha_r$,
and polar part of the Higgs field given at each puncture by
$\sum_1^n B_i \frac{dz}{z^i}$, with the $B_i$ diagonal matrices.
Then there exists a hermitian metric on $\cE$ such that the induced
connection $A=\db^{\cE}+\partial^{\cE}+\theta+\theta^*$ belongs to a
space $\cA$ of connections with datas (\ref{f_alpha}) and
(\ref{f_theta}) at the punctures. The bundle $(\cE,\theta)$ can be
recovered from $(\db_A,\theta_A)$ as its canonical extension.
\end{lemma}
\begin{proof}
The problem consists in constructing an initial metric $h$ on $\cE$.
In order to simplify the ideas, we will restrict to the case
where $B_n$ is regular semisimple, but the general case is similar.
Take a basis $(\sigma_i)$ of
eigenvectors of $\cE_p$ at the puncture $p$, and extend it
holomorphically in a neighborhood. The action of a \emph{holomorphic}
gauge transformation $g=e^u$ of $\cE$ on the Higgs field is by
$$ \theta \longrightarrow g\theta g^{-1}=\exp(\ad u)\theta. $$
From this it is easy to see that by a careful choice of $u$, one can
kill the off-diagonal coefficients of $\theta$ to any finite order.
Therefore we can suppose that in the basis $(\sigma_i)$, we have
$$
\theta = \sum_1^n B_i \frac{dz}{z^i} + \vartheta,
$$
where $\vartheta$ is holomorphic, and the off-diagonal coefficients of
$\vartheta$ vanish up to any fixed order.
Now we choose the flat metric
$$
h=\left(
\begin{array}{ccc}|z|^{2\alpha_1} & & \\ & \ddots & \\ & & |z|^{2\alpha_r}
\end{array}\right).
$$
It is clear that in the orthonormal basis $e_i=s_i/|z|^{\alpha_i}$, we
get exactly the flat model (\ref{f_modDol}), that is
$$
A=\db^{\cE}+\partial^{\cE}+\theta+\theta^*=
d + \re(A_1)id\theta+
\frac{1}{2}\sum_1^n \big(A_i \frac{dz}{z^i} + A_i^* \frac{d\zb}{\zb^i}\big)
  - \beta \frac{dr}{r} + a,
$$
with notations as in section \ref{s_model}.
In the perturbation $a$, the diagonal terms are $C^\infty$, and
the off-diagonal terms can be taken to vanish to any fixed high order.
In particular, we get an element of $\cA$. Actually we have obtained
much more, because $F_A$ vanishes up to any fixed high order (the
diagonal terms do not contribute to the curvature).
\end{proof}

\begin{remark}\label{r_clss}
The extension property for the holomorphic bundle alone (no Higgs
field) follows from earlier work: the connection $A^+$ has curvature
in some $L^p$ for $p>1$ and this implies that a canonical extension as
above exists \cite{Bi92}.
Actually we also need to make precise the singularity of the Higgs field in
the extension: this requires the calculations above.
\end{remark}

Finally, we prove that stability on both sides of (\ref{a_ext}) coincide,
transforming the equivalence of categories (\ref{a_ext}) into an
isomorphism of the moduli spaces, 
$$ \cM_{Dol,an}^s \stackrel{\sim}{\longrightarrow} \cM_{Dol,alg}^s. $$
First we have to introduce an algebraic notion of stability.
A meromorphic Higgs bundle $(\cE,\theta)$ as above has a parabolic degree
$$ \pdeg^{alg}\cE = c_1(\cE)[X] + \sum \alpha_i; $$
if there are several marked points, one must add the contribution of
the weights of each marked point.
Subbundles of $\cE$ also inherit a parabolic degree, and this enables
to define stability.

\begin{lemma}
Let $(\db_A,\theta_A)\in\cA$ be a Higgs bundle and $(\cE,\theta)$ its
canonical extension on $X$. Then $(\db_A,\theta_A)$ is analytically
stable if and only if $(\cE,\theta)$ is algebraically stable.
\end{lemma}
\begin{proof}
The point is to prove that a holomorphic $L^{1,2}$-subbundle (stable
under the Higgs field) extends into an algebraic subbundle of $\cE$,
and that the algebraic and analytic degrees coincide.
Because of remark \ref{r_clss}, this is a consequence of the same
statements in \cite{Si90}.
\end{proof}

\section{The DeRham moduli space}\label{s_DR}

We will not give any detail here, since this is completely parallel to
the results of section \ref{s_Dol}.
We only prove the following technical lemma, which is necessary for
solving the $\db$-problem on components of $\End(E)_k$ with $k\geq 2$.

\begin{lemma}\label{l_db2}
Take $p>2$, $k>1$ and $\delta\in\bR-\bZ$. Then the problem
$$ 
\frac{\partial f}{\partial\zb} + \frac{\lambda}{\zb^k} f = g
$$
has a solution $f=Tg$ such that
$$ \|f\|_{C^0_{\delta}} \leq c \|g\|_{L^p_{-1+\delta}}. $$
The constant $c$ does not depend on $\lambda$.
\end{lemma}
\begin{proof}
The function
$$
\varphi=e^{\frac{\lambda}{(k-1)\zb^{k-1}}-\frac{\bar{\lambda}}{(k-1)z^{k-1}}}
$$
is a solution of
$ \frac{\partial f}{\partial\zb} + \frac{\lambda}{\zb^k} f = 0 $, such
that $|\varphi|=1$. 
Let $T_0$ be the inverse defined by lemma $\ref{l_db1}$, then
$$ Tg = \varphi T_0(\varphi^{-1}g) $$
satisfies the requirements of the lemma.
\end{proof}

An element of the moduli space $\cM_{DR,an}$ is represented by a flat
connection $A\in\cA$. The holomorphic bundle underlying $A$ has the
$\db$-operator $D_A^{0,1}$. As in the previous section, relying on lemmas
\ref{l_db1} and \ref{l_db2}, near a puncture it is possible to produce
a complex gauge transformation $g$ such that 
$$ g(D_A^{0,1}) = D_0^{0,1} = \db +
\frac{1}{2} \sum_2^n A_i^* \frac{d\zb}{\zb^i} -
 \frac{\beta+\im A_1}{2} \frac{d\zb}{\zb} .
$$
Therefore we have $D_A^{0,1}$-holomorphic sections $(\tau_i)$ given by
$$
\tau_i = |z|^{\beta_i+\im\mu_i} g^{-1} 
         \exp\bigg(
  \sum_2^n \frac{A_i^*}{2(i-1)\zb^{i-1}} - \frac{A_i}{2(i-1)z^{i-1}}
             \bigg) e_i,
$$
and this basis of holomorphic sections define a canonical holomorphic
extension $\cF$ over the puncture.

With respect to this extension, $A$ becomes an integrable connection
with irregular singularities, and the polar part of $A$ remains equal to
that of the model,
$$ d + A_n \frac{dz}{z^n} + \cdots + A_1 \frac{dz}{z}, $$
so that we get an element of the moduli space $\cM_{DR,alg}$ of such
integrable connections.

Also the extension has a parabolic structure with weights $\beta_i$,
and this enables us to define a parabolic degree $\pdeg^{alg}\cF$ and
therefore the algebraic stability of $(\cF,A)$.

We finally have all the ingredients of the isomorphism
$$ \cM_{DR,an}^s \longrightarrow \cM_{DR,alg}^s. $$

\begin{remark}\label{rem:weights}
The weights $\gamma_i$ of the local system are the order of growth of
parallel sections on rays going to the singularity. From the above formula,
they are equal to $$\gamma_i=\beta_i-\re\mu_i=-2\re\lambda_i.$$
By \cite[proposition 11.1]{Bi97}, the parabolic degree of $\cF$ is 
$$ \pdeg^{alg}\cF=\sum\gamma_i, $$
where the sum has to be understood for all punctures.
In particular, if all weights $\gamma_i$ are taken to be zero, then
the same is true for subbundles, so the degree for subbundles is
always zero, so that stability reduces to irreducility of the connection.
\end{remark}

\subsubsection*{Sufficient stability conditions}

We will describe some simple conditions on the
parameters such that all points of $\cM_{DR}$ are stable.
Suppose $A$ is a meromorphic connection on a holomorphic 
vector bundle $\cE\to X$ as constructed from the extension
procedure
above.
Thus in some local trivialisation near the $i$th singularity
the polar part of
$A$ takes the form of the model
$$ d + \iA_{n_i} \frac{dz_i}{z_i^{n_i}} + \cdots + \iA_1
\frac{dz_i}{z_i}, $$
where $z_i$ is a local coordinate and $\iA_j$ are diagonal
matrices.
We wish to assume now that all of the leading coefficients
$\iA_{n_i}$ are
regular (have distinct eigenvalues).
Note that the eigenvalues of the residues
$\iA_1$ are uniquely determined by $A$ upto order
(independent of the
coordinate choice). 
Now if $\cF$ is a subbundle of $\cE$ preserved by $A$
we may choose a trivialisation of $\cF$ putting the induced
connection
on $\cF$ in model form (with residues $\iB_1$ say)
and then extend this to a trivialisation of $\cE$ as above.
In particular it follows that the eigenvalues of $\iB_1$ are
a subset
of the eigenvalues of $\iA_1$.
However (by considering the induced connection on
$\det(\cE)$) we know 
that the (usual) degree of $\cE$ is minus the sum of traces of
the residues:
$$\deg(\cE) = - \sum_i\tr(\iA_1)$$
and similarly for $\cF$.
Thus we can ensure that $A$ has no proper nontrivial
subconnections
by choosing the models for $A$ such that none of the (finite
number of) ``subsums'' 
\begin{equation}\label{subsums}
\sum_i\sum_{j\in S_i} (\iA_1)_{jj}
\end{equation}
of the residues 
is an integer, where $S_i\subset\{1,\ldots,\rank(\cE)\}$ are
finite subsets
of size $k$ and $k$ ranges from $1$ to $\rank(\cE)-1$. 
Thus under such (generic) conditions any such connection $A$
is stable.

\subsubsection*{An example}

Let $\cM$ be a moduli space of integrable connections with
two poles 
on the projective
line $P^1$ of order two at zero and order one at infinity.
Consider the subspace $\cM^*\subset\cM$ of connections such
that the
underlying holomorphic bundle is trivial. Therefore points
of $\cM^*$
are (globally) represented by connections of the form
\begin{equation}\label{eqn: global connection}
d + A_0 \frac{dz}{z^2} + B \frac{dz}{z}. 
\end{equation}
We will assume that $A_0$ is diagonal with distinct
eigenvalues and
that none of the eigenvalues of $B$ differ by integers but allow arbitrary
parabolic weights.
It then follows that the model connections at zero and
infinity are
\begin{equation}\label{models}
d+ A_0 \frac{dz}{z^2} + \Lambda \frac{dz}{z}, \qquad d+B_0\frac{dz}{z}
\end{equation}
respectively, where $\Lambda$ is the diagonal part of $B$
and $B_0$ is
a diagonalisation of $B$.
Once these models are fixed we see that $B$ is restricted to
the
adjoint orbit $\cO$ containing $B_0$ and has diagonal part
fixed to
equal $\Lambda$. 
Since $A_0$ is fixed and regular the 
remaining gauge freedom in \eqref{eqn: global connection} is
just
conjugation by $T^\bC$.
Now if we use the trace to identify $\cO$ with a coadjoint
orbit and so
give it a complex symplectic structure, then the map 
$\delta:\cO\to\lt_\bC$ taking $B$ to its diagonal part is a
moment map
for this torus action and so we have an isomorphism
$$\cM^*\cong \cO\spq T^\bC$$
of the moduli space with the complex symplectic quotient at
the value
$\Lambda$ of the moment map.

On the other hand, the same symplectic quotient underlies a
complete
hyperK{\"a}hler metric obtained by taking the 
hyperK{\"a}hler quotient of Kronheimer's
hyperK{\"a}hler metric \cite{Kr90} on $\cO$ by the maximal
compact torus. 
Nevertheless, in general this quotient metric
on $\cM^*$ does not coincide with the metric of theorem
\ref{t_m},
because that is a complete metric on $\cM$, which is larger,
as we
will show below. 
Therefore, varying $A_0$ in the regular part of the Cartan
subalgebra leads to a family of hyperK{\"a}hler metrics
on $\cO\spq T^\bC$ which become complete in a larger space.
We remark that in this example the full space $\cM$ may be
analytically
identified (cf. \cite{smpgafm,saqh}) with the complex
symplectic
quotient $\cL\spq T^\bC$ of a symplectic leaf $\cL\subset
G^*$ of the
simply connected Poisson Lie group $G^*$ dual to
$GL_r(\bC)$.

\begin{lemma}\label{lem: nontrivial bundles}
There are stable connections on non-trivial bundles with
models of the type
considered in the above example.
\end{lemma}
\begin{proof}
We will do this in the rank three case (this is the simplest
case
since then $\dim_\bC \cM=2$; one may easily generalise to
higher rank).
Suppose we have 
$g\in GL_3(\bC)$ and diagonal matrices $A_0,B'_0$ such that
$A_0,e^{2\pi i B'_0}$ have distinct eigenvalues and
\begin{enumerate}
\item[1)] the  matrix entry $(g^{-1}A_0g)_{31}$ is zero, and
\item[2)] the pair of diagonal matrices $-B'_0,
\Lambda:=\delta(gB'_0g^{-1})$ have no
integral subsums (in the sense of \eqref{subsums}).
\end{enumerate}
Then consider the meromorphic connection on the
bundle 
$\cO(1)\oplus\cO\oplus\cO(-1)\to P^1$ defined by the
clutching map
$h=\diag(z,1,z^{-1})$ and equal to 
$$d+ g^{-1}A_0g \frac{dz}{z^2} + B'_0 \frac{dz}{z}$$
on $\bC\subset P^1$.
Assumption 1) implies this is equivalent to the models
\eqref{models}
at $0$ and $\infty$, with $B_0= B'_0+ \diag(1,0,-1)$.
Then assumption 2) implies it is stable.
Finally one may easily construct such matrices by observing
that
companion matrices have zero $(31)$ entry. For example
$$A_0= \lambda\left(
\begin{matrix}
1 &  &  \\
 & 2 &  \\
 &  & 4 
\end{matrix}
\right),\qquad 
B_0=\left(
\begin{matrix}
p &  &  \\
 & q &  \\
 &  & r 
\end{matrix}
\right), \qquad
g=\left(
\begin{matrix}
1 & 1 & 1 \\
2 & 4 & 8 \\
3 & 12 & 48 
\end{matrix}
\right)$$
have the desired properties if $\lambda$ is nonzero
and $p,q,r$
are sufficiently generic (e.g. if $\{1,p,q,r\}$ are linearly
independent
over the rational numbers).
\end{proof}

\section{Proof of theorem \ref{t_he}}\label{s_proof}

We will prove only the most difficult way, that is the isomorphism with the
Dolbeault moduli space.
We start with a stable Higgs bundle $(\db_{A_0},\phi_{A_0})\in\cA$,
satisfying the integrability condition $(D_{A_0}'')^2=0$, that is
$$ \db_{A_0}\theta_{A_0}=0, $$
and we try to find a complex gauge transformation $g\in\cG_{\bC}$ taking
it to a solution $A=g(A_0)$ of the selfduality equations, that is
satisfying the additional equation
$$ F_A=0. $$
This problem is equivalent to finding the Hermitian-Einstein metric
$h=g^*g$ on the Higgs bundle.

The idea, as in \cite[section 8]{Bi97} is to minimize the Donaldson
functional $M(h)$ for $h=g^*g$ under the constraint
\begin{equation}\label{constraint}
\|\Lambda F_{(\db_{A_0},\phi_{A_0})}^h\|_{L^2_{-2+\delta}} \leq B.
\end{equation}

We need the two following technical lemmas.
\begin{lemma}\label{l_t1}
If $f\in L^{2,2}_{-2+\delta}$ and $g\in L^2_{-2+\delta}$ are positive
functions such that $\Delta f\leq g$ and $\|g\|_{L^2_{-2+\delta}}\leq B$, then
$$ \|f\|_{C^0} \leq c_0(B) + c_1(B) \|f\|_{L^1}, $$
and $\|f\|_{C^0}$ goes to $0$ when $\|f\|_{L^1}$ goes to zero.
\end{lemma}
\begin{proof}
This follows easily (see \cite[proposition 2.1]{Si88}) from the fact
that if $v\in L^2_{-2+\delta}$, then the problem $\Delta u=v$ on the
unit disk, with Dirichlet boundary condition, can be solved with $u\in
L^{2,2}_{-2+\delta}\subset C^0$. 
\end{proof}

\begin{lemma}\label{l_t2}
If we have a sequence of metrics $h_j=h_0u_j$ with
\begin{enumerate}
\item $h_j$ has a $C^0$-limit $h_\infty$;
\item $u_j\in L^{2,2}_{-2+\delta}$ and $h_j$ satisfies the constraint
(\ref{constraint}); 
\item $\|D''_{A_0}u_j\|_{L^2}$ is bounded;
\end{enumerate}
then the limit is actually a $L^{2,2}_{-2+\delta}$-limit:
$h_\infty=h_0 u_\infty$ and $u_\infty\in L^{2,2}_{-2+\delta}$.
\end{lemma}
\begin{proof}
This is a local statement. Outside the puncture, the statement is
proven for example in \cite[lemma 6.4]{Si88}, the point here is to
prove that $C^0$-convergence implies the convergence in our weighted
Sobolev spaces.

Note $D=D_{A_0}$, $D''=D''_{A_0}$, etc. We use the freedom, from the
proof of lemma \ref{l_initial}, to choose an initial metric $h$ with
$\Lambda F^h$ bounded.
Observe that because of the formula
$$ (D'')^*D''-(D')^*D'= i\Lambda F^h, $$
the hypothesis on $\|D''u_j\|_{L^2}$ also implies that
$\|D'u_j\|_{L^2}$ is bounded.

We have the formula
$$
(D')^*D'u_j = iu_j(\Lambda F^{h_j}-\Lambda F^h)
             +i\Lambda (D''u_j)u_j^{-1}(D'u_j).
$$

Let $\chi$ be a cut-off function, with compact support in the disk,
such that $\chi|_{\Delta_{1/2}}=1$. Then we obtain
\begin{multline}\label{ddchiuj}
(D')^* D' \chi u_j =
  i\chi u_j(\Lambda F^{h_j}-\Lambda F^h)
  +i\Lambda (D''u_j)u_j^{-1}(D'\chi u_j)\\
  +d\chi \odot Du_j + \frac{1}{2}(\Delta\chi) u_j
\end{multline}
In particular, we get
\begin{equation}\label{du}
\|(D')^* D' \chi u_j\|_{L^2_{-2+\delta}}
\leq c \big(1 + \|D''u_j\|_{L^4_{-1}} \|D'\chi u_j\|_{L^4_{-1+\delta}} \big)
\end{equation}
Notice that, by corollary \ref{c_control}, because $(D')^2=0$, we have
$$
\|(D')^* D'\chi u_j\|_{L^2_{-2+\delta}}
\geq c\|D'\chi u_j\|_{L^{1,2}_{-2+\delta}}
\geq c\|D'\chi u_j\|_{L^4_{-1+\delta}}
$$
and therefore from (\ref{du})
\begin{equation}\label{chiuj}
\|D'\chi u_j\|_{L^{1,2}_{-2+\delta}} (1-c\|D''u_j\|_{L^4_{-1}})
\leq c
\end{equation}

The important point here is that all constants are invariant by
homothety. Indeed the $L^2$-norm of 1-forms is conformally invariant,
so $\|(D''\oplus D')u\|_{L^2}$ remains bounded; the same is true for the
$L^4_{-1}$-norm of 1-forms; in corollary \ref{c_control}, the
constants do not depend on an homothety; finally, in the Sobolev embedding
$$
\|d(\chi f)\|_{L^2_{-2+\delta}}
\geq c \|\chi f\|_{L^4_{-1+\delta}},
$$
the norms on both sides are rescaled by the same factor under an homothety.

Now suppose that there exists some disk $\Delta_\rho$ of radius $\rho$
such that for all $j$ one has
$$
\|D''u_j\|_{L^4_{-1}(\Delta_\rho)} < \frac{1}{2c} ;
$$
then, because this norm is invariant under homothety, we can rescale
the disk $\Delta_\rho$ by the homothety $h_{\rho}$ to the unit disk
$\Delta$, and applying (\ref{chiuj}) we get that
$\chi h_{\rho}^*u_j$ is bounded in $L^{2,2}_{-2+\delta}$, and we get
the lemma.

Now suppose on the contrary that there exist radii $\rho_j\rightarrow 0$
such that
$$
\|D''u_j\|_{L^4_{-1}(\Delta_{\rho_j})} \geq \frac{1}{2c} ;
$$
by taking some smaller $\rho_j$, one can arrange it so that actually
\begin{equation}\label{arrange}
\|D''u_j\|_{L^4_{-1}(\Delta_{\rho_j})}
+\|D'u_j\|_{L^4_{-1}(\Delta_{\rho_j})}
= \frac{1}{2c}.
\end{equation}

We will see that this hypothesis leads to a contradiction.
First we prove that the ``energy'' in (\ref{arrange})
cannot concentrate near the origin.
Again using the homothety $h_{\rho_j}$, we deduce from (\ref{chiuj})
that $D'(\chi h_{\rho_j}^*u_j)$ is bounded in $L^{1,2}_{-2+\delta}$;
of course one must be a bit careful here, because these
$L^{1,2}$-norms do depend on $j$, so we actually mean (write
$D_j=h_{\rho_j}^* D$)
\begin{equation}\label{e_bdd}
 \|D'_j(\chi h_{\rho_j}^*u_j)\|_{L^2_{-1+\delta}}
  +\|(\nabla^{D^+_j}\oplus \phi_j)D'_j(\chi
h_{\rho_j}^*u_j)\|_{L^2_{-2+\delta}} \leq C ;
\end{equation}
in particular,
$$ \|d\big|D'_j(\chi h_{\rho_j}^*u_j)\big|\|_{L^2_{-2+\delta}} $$
is bounded, and because of the compact inclusion
$L^{1,2}_{-2+\delta}\subset L^4_{-1}$ (see remark \ref{r_cpct}),
we deduce that the functions $|D'_j(\chi h_{\rho_j}^*u_j)|$
converge strongly in $L^4_{-1}$ to a limit (and the same is true for
$|D''_j(\chi h_{\rho_j}^*u_j)|$).

Actually, we can deduce a bit more from (\ref{e_bdd}): indeed, the
operators $D_j=h_{\rho_j}^* D$ become very close to the model
$h_{\rho_j}^* D_0$ when $j$ goes to infinity, so it is enough to suppose
that $D=D_0$ near the puncture; now for components $u(k)$ for $k\geq
2$, we have by (\ref{e_bdd})
$$
\frac{\lambda_k}{\rho_j^{k-1}}
\| D'_j(\chi h_{\rho_j}^*u_j(k))\|_{L^2_{-1+\delta}}
\leq
\| \phi_j\otimes D'_j(\chi h_{\rho_j}^*u_j(k)) \|_{L^2_{-2+\delta}}
\leq C ;
$$
and therefore the limit in $L^4_{-1}$ of
$|D'_j(\chi h_{\rho_j}^*u_j(k))|$ must be zero.

Therefore we are left with only the limit of the components with $k=0$
or $k=1$: observe now that on these components, the operator $D_0$ is
homothety invariant, so it makes sense to look at the limit of the
operators $D_j$.
Recall that $u_j$ has a $C^0$-limit, and $\rho_j\rightarrow 0$, so
that $h_{\rho_j}^*u_j$ has a constant limit.
We deduce that actually, for components of $u_j$ with $k=0$ or $1$,
the limit  in $L^{4}_{-1}$ of $(D''_j\oplus D'_j)(h_{\rho_j}^*\chi u_j)$
must be $0$, which implies
$$
\|D''u_j\|_{L^4_{-1}(\Delta_{\frac{1}{2}\rho_j})}
+\|D'u_j\|_{L^4_{-1}(\Delta_{\frac{1}{2}\rho_j})}
 \rightarrow 0.
$$

Therefore we have proven, as announced, that the ``energy''
(\ref{arrange}) cannot concentrate near the origin.
We deduce that there exist points
$$x_j\in \Delta_{\rho_j}-\Delta_{\frac{1}{2}\rho_j},$$
such that
\begin{equation}\label{lbdu}
\|D''u_j\|_{L^4_{-1}(\Delta_{\frac{1}{8}\rho_j}(x_j))} \geq \frac{1}{100c}.
\end{equation}
Now we rescale the disk $\Delta_{\frac{1}{4}\rho_j}(x_j)$ centered at
$x_j$ into the unit disk $\Delta$ by a homothety $h_j'$;
a similar argument gives from (\ref{ddchiuj}) the estimate
$$
\|\chi (h'_j)^*u_j\|_{L^{2,2}}
  \leq c (1+\|\nabla (h'_j)^*u_j\|_{L^4}^2)
  \leq C.
$$
Hence we can extract a strongly $L^{1,4}$ convergent subsequence 
$(h'_j)^*u_j$, but the limit must again be flat, and this contradicts
the fact that by (\ref{lbdu}) the $L^4$-norm of
$(D'\oplus D'')\big((h'_j)^*u_j\big)$ on $\Delta_{\frac{1}{2}}$ is
bounded below.
\end{proof}

Using these two lemmas, the proof is now a standard adaptation of that
in \cite{Bi97}. Indeed, from lemma \ref{l_t1} one deduces by Simpson's
method that if the bundle is stable, then Donaldson's functional is
bounded below, and a minimizing sequence $h_j=hu_j$ must converge in $C^0$ to
some limit; moreover, $\|D''u_j\|_{L^2}$ remains bounded. Lemma
\ref{l_t2} gives the stronger convergence in Sobolev spaces
$L^{2,2}_{-2+\delta}$, and one can then deduce that the limit actually
solves the equation.\qed

\hfill \textsc{Irma}, Universit{\'e} Louis Pasteur et \textsc{Cnrs}

\hfill 7 rue Ren{\'e} Descartes, 67084 \textsc{Strasbourg Cedex, France}


\begin{thebibliography}{Boa01b}

\bibitem[Biq91]{Bi91}
O.~Biquard.
\newblock Fibr\'es paraboliques stables et connexions singuli\`eres plates.
\newblock {\em Bull. Soc. Math. France}, 119(2):231--257, 1991.

\bibitem[Biq92]{Bi92}
O.~Biquard.
\newblock Prolongement d'un fibr\'e holomorphe hermitien \`a\ courbure ${L}\sp
  p$ sur une courbe ouverte.
\newblock {\em Internat. J. Math.}, 3(4):441--453, 1992.

\bibitem[Biq97]{Bi97}
O.~Biquard.
\newblock Fibr\'es de {H}iggs et connexions int\'egrables: le cas logarithmique
  (diviseur lisse).
\newblock {\em Ann. Sci. \'Ecole Norm. Sup. (4)}, 30(1):41--96, 1997.

\bibitem[Boa]{saqh}
P.~P. Boalch.
\newblock Quasi-{H}amiltonian geometry of meromorphic connections.
\newblock \texttt{math.DG/0203161}.

\bibitem[Boa01a]{smpgafm}
P.~P. Boalch.
\newblock {S}tokes matrices, {P}oisson {L}ie groups and {F}robenius manifolds.
\newblock {\em Invent. math.}, 146(3):479--506, 2001.

\bibitem[Boa01b]{Bo99}
P.~P. Boalch.
\newblock Symplectic manifolds and isomonodromic deformations.
\newblock {\em Adv. Math.}, 163(2):137--205, 2001.

\bibitem[Bot95]{Bottacin}
F.~Bottacin.
\newblock Symplectic geometry on moduli spaces of stable pairs.
\newblock {\em Ann. Sci. Ecole Norm. Sup. (4)}, 28(4):391--433, 1995.

\bibitem[CK01]{CheKap01}
S.~Cherkis and A.~Kapustin.
\newblock Nahm transform for periodic monopoles and $\mathscr{N}=2$ super
  {Y}ang-{M}ills theory.
\newblock {\em Comm. Math. Phys.}, 218(2):333--371, 2001.

\bibitem[CK02]{CheKap}
S.~Cherkis and A.~Kapustin.
\newblock Hyper-{K}\"ahler metrics from periodic monopoles.
\newblock {\em Phys. Rev. D (3)}, 65(8):084015 (10 p.), 2002.

\bibitem[Don87]{Do87}
S.~K. Donaldson.
\newblock Twisted harmonic maps and the self-duality equations.
\newblock {\em Proc. London Math. Soc. (3)}, 55(1):127--131, 1987.

\bibitem[Hit87]{Hi87}
N.~J. Hitchin.
\newblock The self-duality equations on a {R}iemann surface.
\newblock {\em Proc. London Math. Soc. (3)}, 55(1):59--126, 1987.

\bibitem[Kro90]{Kr90}
P.~B. Kronheimer.
\newblock {A hyper-K{\"a}hlerian structure on coadjoint orbits of a semisimple
  complex group.}
\newblock {\em J. Lond. Math. Soc., II. Ser.}, 42(2):193--208, 1990.

\bibitem[LM85]{LoMO85}
R.~B. Lockhart and R.~C. McOwen.
\newblock Elliptic differential operators on noncompact manifolds.
\newblock {\em Ann. Scuola Norm. Sup. Pisa Cl. Sci. (4)}, 12(3):409--447, 1985.

\bibitem[Mar94]{Markman}
E.~Markman.
\newblock Spectral curves and integrable systems.
\newblock {\em Comp. Math.}, 93(3):255--290, 1994.

\bibitem[MR91]{MarRam91}
J.~Martinet and J.-P. Ramis.
\newblock Elementary acceleration and multisummability. {I}.
\newblock {\em Ann. Inst. H. Poincar\'e Phys. Th\'eor.}, 54(4):331--401, 1991.

\bibitem[Sab99]{Sa99}
C.~Sabbah.
\newblock Harmonic metrics and connections with irregular singularities.
\newblock {\em Ann. Inst. Fourier (Grenoble)}, 49(4):1265--1291, 1999.

\bibitem[Sim88]{Si88}
C.~T. Simpson.
\newblock Constructing variations of {H}odge structure using {Y}ang-{M}ills
  theory and applications to uniformization.
\newblock {\em J. Amer. Math. Soc.}, 1(4):867--918, 1988.

\bibitem[Sim90]{Si90}
C.~T. Simpson.
\newblock Harmonic bundles on noncompact curves.
\newblock {\em J. Amer. Math. Soc.}, 3(3):713--770, 1990.

\bibitem[Sim92]{Si92}
C.~T. Simpson.
\newblock Higgs bundles and local systems.
\newblock {\em Inst. Hautes \'Etudes Sci. Publ. Math.}, 75:5--95, 1992.

\end{thebibliography}
\end{document}